\documentclass[11pt]{amsart}

\usepackage{latexsym}
\usepackage{amssymb}
\usepackage{amscd}
\usepackage[all]{xy}
\usepackage[mathscr]{euscript}
\usepackage{dsfont}
\usepackage{hyperref}
\usepackage{graphicx}
\usepackage{amsfonts}
\usepackage{amsmath}
\usepackage{amsthm}
\usepackage{aliascnt}

\normalsize

\RequirePackage{xspace}
\newcommand{\thought}[1]{}
\renewcommand{\thought}[1]{ \textbf{[#1]}}

\usepackage{enumerate}        % Fancy enumerations;  used for Roman numbering.
\newenvironment{roenumerate}{\begin{enumerate}[\upshape (i)]}{\end{enumerate}}
\newcommand\nc {\newcommand}
\newcommand\rnc{\renewcommand}

\newcount\blopone
\newcount\xone
\newcount\xtwo
\newcount\ytwo

\newtheorem{theorem}{Theorem}[section]
\newtheorem{prop}[theorem]{Proposition}
\newtheorem{com}[theorem]{Comment}
\newtheorem{apl}[theorem]{Application}
\newtheorem{exercise}[theorem]{Exercise}
\newtheorem{redu}[theorem]{Reduction}
\newtheorem{refinement}[theorem]{Refinement}
\newtheorem{summary}[theorem]{Summary}
\newtheorem{importnota}[theorem]{Important Notation}
\newtheorem{prblm}[theorem]{Problem}
\newtheorem{notation}[theorem]{Notation}
\newtheorem{explanation}[theorem]{Explanation}
\newtheorem{defin}[theorem]{Definition}
\newtheorem{caution}[theorem]{Caution}
\newtheorem{remark}[theorem]{Remark}
\newtheorem{reminder}[theorem]{Reminder}
\newtheorem{illustration}[theorem]{Illustration}
\newtheorem{observation}[theorem]{Observation}
\newtheorem{lemma}[theorem]{Lemma}
\newtheorem{construction}[theorem]{Construction}
\newtheorem{discussion}[theorem]{Discussion}
\newtheorem{corollary}[theorem]{Corollary}
\newtheorem{example}[theorem]{Example}
\newtheorem{conclusion}[theorem]{Conclusion}
\newtheorem{sketch}[theorem]{Sketch}
\newtheorem{triviality}[theorem]{Triviality}
\newtheorem{proto}[theorem]{Prototype Quasifibration}
\newtheorem{cauex}[theorem]{Cautionary Example}
\newtheorem{hypo}[theorem]{Hypothesis}
\newtheorem{subth}{ }[theorem]
\newtheorem{case}{Case}[theorem]
\newtheorem{ssubth}{ }[subth]
\newtheorem{facts}[theorem]{Facts}
\newtheorem{history}[theorem]{Historical Survey}
\newtheorem{proofs}[theorem]{Discussion of the Proofs, Old and New}
\newtheorem{discl}[theorem]{Disclaimer}

\nc\tri[1]{\begin{triviality}
\label{#1}}
\nc\fac[1]{\begin{facts}
\label{#1}
\begin{em}}
\nc\app[1]{\begin{apl}
\label{#1}
\begin{em}}
\nc\skt[1]{\begin{sketch}
\label{#1}
\begin{em}}
\nc\hst[1]{\begin{history}
\label{#1}
\begin{em}}
\nc\pfs[1]{\begin{proofs}
\label{#1}
\begin{em}}
\nc\cas[1]{\begin{case}
\label{#1}
\begin{em}}
\nc\rfn[1]{\begin{refinement}
\label{#1}}
\nc\prt[1]{\begin{proto}
\label{#1}}
\nc\lem[1]{\begin{lemma}
\label{#1}}
\nc\pro[1]{\begin{prop}
\label{#1}}
\nc\thm[1]{\begin{theorem}
\label{#1}}
\nc\dis[1]{\begin{discussion}
\label{#1}
\begin{em}}
\nc\dsc[1]{\begin{discl}
\label{#1}
\begin{em}}
\nc\cor[1]{\begin{corollary}
\label{#1}}
\nc\dfn[1]{\begin{defin}
\label{#1}}
\nc\sthm[1]{\begin{subth}
\label{#1}}
\nc\exm[1]{\begin{example}
\label{#1}
\begin{em}}
\nc\obs[1]{\begin{observation}
\label{#1}
\begin{em}}
\nc\plm[1]{\begin{prblm}
\label{#1}
\begin{em}}
\nc\rmk[1]{\begin{remark}
\label{#1}
\begin{em}}
\nc\rmd[1]{\begin{reminder}
\label{#1}
\begin{em}}
\nc\ntn[1]{\begin{notation}
\label{#1}
\begin{em}}
\nc\exe[1]{\begin{exercise}
\label{#1}
\begin{em}}
\nc\xpl[1]{\begin{explanation}
\label{#1}
\begin{em}}
\nc\smr[1]{\begin{summary}
\label{#1}
\begin{em}}
\nc\cau[1]{\begin{caution}
\label{#1}
\begin{em}}
\nc\hyp[1]{\begin{hypo}
\label{#1}}
\nc\imn[1]{\begin{importnota}
\label{#1}
\begin{em}}
\nc\rdn[1]{\begin{redu}
\label{#1}
\begin{em}}
\nc\cax[1]{\begin{cauex}
\label{#1}
\begin{em}}
\nc\cmt[1]{\begin{com}
\label{#1}
\begin{em}}
\nc\con[1]{\begin{construction}
\label{#1}
\begin{em}}
\nc\ill[1]{\begin{illustration}
\label{#1}
\begin{em}}
\nc\ssthm[1]{\begin{ssubth}
\label{#1}
\begin{em}}
\nc\cnc[1]{\begin{conclusion}
\label{#1}
\begin{em}}

\nc\elem{\end{lemma}}
\nc\erdn{\end{em}\end{redu}}
\nc\erfn{\end{refinement}}
\nc\eprt{\end{proto}}
\nc\ethm{\end{theorem}}
\nc\ecor{\end{corollary}}
\nc\edfn{\end{defin}}
\nc\esthm{\end{subth}}
\nc\epro{\end{prop}}
\nc\etri{\end{triviality}}
\nc\eexm{\end{em}
\end{example}}
\nc\eobs{\end{em}
\end{observation}}
\nc\ecmt{\end{em}
\end{com}}
\nc\efac{\end{em}
\end{facts}}
\nc\eapp{\end{em}
\end{apl}}
\nc\ermk{\end{em}
\end{remark}}
\nc\ermd{\end{em}
\end{reminder}}
\nc\eill{\end{em}
\end{illustration}}
\nc\eplm{\end{em}
\end{prblm}}
\nc\ecas{\end{em}
\end{case}}
\nc\eskt{\end{em}
\end{sketch}}
\nc\ecau{\end{em}
\end{caution}}
\nc\ecax{\end{em}
\end{cauex}}
\nc\eimn{\end{em}
\end{importnota}}
\nc\entn{\end{em}
\end{notation}}
\nc\eexe{\end{em}
\end{exercise}}
\nc\expl{\end{em}
\end{explanation}}
\nc\edis{\end{em}
\end{discussion}}
\nc\edsc{\end{em}
\end{discl}}
\nc\econ{\end{em}
\end{construction}}
\nc\esmr{\end{em}
\end{summary}}
\nc\ehst{\end{em}
\end{history}}
\nc\epfs{\end{em}
\end{proofs}}
\nc\ehyp{
\end{hypo}}
\nc\ecnc{\end{em}
\end{conclusion}}
\nc\essthm{\end{em}
\end{ssubth}}

\nc\sst{\scriptstyle}
\newcommand{\comment}[1]{}
\newcommand{\ri}{\longrightarrow}

\newcommand{\zz}{{\mathbb Z}}

\newcommand{\nn}{{\mathbb N}}

\newcommand{\D}{{\mathbf D}}

\nc\op{^{\hbox{\rm\tiny op}}}
\nc\mth{^{\hbox{\rm\tiny th}}}

\nc\script{\mathscr}
\nc\z{\zeta}
\nc\bc{{\mathbb{BC}}}
\nc\ct{{\script T}}
\nc\cf{{\script F}}
\nc\cg{{\script G}}
\nc\ch{{\script H}}
\nc\ck{{\script K}}
\nc\cl{{\script L}}
\nc\cv{{\script V}}
\nc\ce{{\script E}}
\nc\cm{{\script M}}
\nc\cn{{\script N}}
\nc\cs{{\script S}}
\nc\car{{\script R}}
\nc\cd{{\script D}}
\nc\cc{{\script C}}
\nc\ca{{\script A}}
\nc\ci{{\script I}}
\nc\co{{\script O}}
\nc\cu{{\script U}}
\nc\cx{{\script X}}
\nc\cy{{\script Y}}
\nc\cz{{\script Z}}
\nc\Cp{{\script P}}
\nc\cq{{\script Q}}
\nc\bd{\begin{description}}
\nc\ed{\end{description}}
\nc\ctob{{\script C}at\big(\ci^{op},\ca\big)}
\nc\clim{{\ds\mathop{\rm lim}_{\ds\longleftarrow}}\,}
\nc\climone{{\ds{\mathop{\rm lim}_{\ds\longleftarrow}}}^1\,}
\nc\climi{\clim_{\!i}\,}
\nc\climn{\clim^{\!n}\,}
\nc\colim{{\ds\mathop{\rm colim}_{\ds\la}}}
\nc\colimj{{\ds\mathop{\rm colim}_{\ds\la}}{}_{j\,}}
\nc\oa{\overline{\ca}}
\nc\s{\sigma}
\nc\ta{\tau}
\nc\os{\overline\sigma}
\nc\ot{\overline\tau}
\nc\T{\Sigma}
\nc\Tm{\Sigma^{-1}}
\nc\de[1]{{\mathop{\rm deg(#1)}}}
\nc\Ad[1]{\mathop{\rm Ad}(#1)}
\nc\ad[1]{\mathop{\rm ad}(#1)}
\nc\kth{{\it K}--theory}
\nc\loc[1]{{\text{\rm Loc}(#1)}}
\nc\coloc[1]{{\text{\rm Coloc}(#1)}}

\def\der #1 {D\left(#1\right)}
\nc\prf{\begin{proof}}
\nc\eprf{\end{proof}}
\nc\ds{\displaystyle}
\nc\Tor{\text{\rm Tor}}

\nc\cb{{\script B}}
\nc\ab{{\script A}b}

\nc\be{\begin{roenumerate}}
\nc\ee{\end{roenumerate}}

\nc\cat[1]{{\script C}at\Big({\big\{#1\big\}}\op\,\,,\,\,\ab\Big)}
\nc\csab{{\script C}at\big(\cs^{op},\ab\big)}
\nc\ctab{{\script C}at\Big({\{\ct^\alpha\}}^{op},\ab\Big)}
\nc\csex{{\script E}x\big(\cs^{op},\ab\big)}
\nc\ctex{{\script E}x\Big({\{\ct^\alpha\}}^{op},\ab\Big)}
\nc\sub{\qquad\subset\qquad}
\nc\ctr[1]{{\left.\ct\left(-,#1\right)\right|}_{\cs}}
\nc\ctrf[2]{{\left.\ct\left(#1,#2\right)\right|}_{\cs}}
\nc\Ctr[1]{{\left.\ct\left(-,#1\right)\right|}_{\ct^\alpha}}
\nc\Ctrf[2]{{\left.\ct\left(#1,#2\right)\right|}_{\ct^\alpha}}

\nc\la{\longrightarrow}
\nc\nin{\noindent}
\nc\cad[1]{\text{card}(#1)}
\nc\eq{\quad=\quad}
\nc\BA{\begin{array}{c}}
\nc\EA{\end{array}}
\nc\barr{
\[
\begin{array}{cccccccccccccccc}
}
\nc\earr{
\end{array}
\]
}
\nc\as[1]{{\langle S\rangle}^{#1}}
\nc\sh{\text{\it shift}}

\nc\yy[1]{{\left.\ct\left(-,#1\right)\right|}_{\ct^c}}
\nc\vrep[2]{{\left.\ct\left(#1,#2\right)\right|}_{\ct^\alpha}}
\nc\da{\downarrow}
\nc\Hom{{\mathop{\rm Hom}}}
\nc\HHom{{\script H}{\mathop{\rm om}}}
\nc\End{{\mathop{\rm End}}}
\nc\Ext{{\mathop{\rm Ext}}}
\nc\mr\Modtc
\nc\PExt{{\mathop{\rm PExt}}}
\nc\stm{\text{\rm stmod}(kG)}
\nc\stM{\text{\rm StMod}(kG)}
\nc\e{\varepsilon}
\nc\p{\varphi}

\nc\rs{\s^{-1}A}
\nc\br{{\{\s^{-1}A\}}}
\nc\ra\ri
\nc\y[1]{\mathbf{y}#1}
\nc\x[1]{\mathbf{z}#1}
\nc\mmod[1]{#1\text{--\rm mod}}
\nc\Mod[1]{#1\text{--\rm Mod}}
\nc\Md {\ensuremath{\mathop{\textup{Mod}}}}
\rnc\mod[1]{\ensuremath{\mathop{#1\textup{--mod}}}\xspace}
\nc\MMod[1]{\text{\rm Mod-}#1}
\nc\Modtc{\Mod{\ct^c}}
\nc\pgldim[1]{\mathop{\rm pgldim}\,#1}
\nc\tf{{\rm [TR5]}}
\nc\tfs{{\rm [TR5$^*$]}}
\nc\Fun{\text{\rm Funct}(F\op,\ab)}
\nc\sym{\text{\rm Sym}}
\nc\sgn{\text{\rm sgn}}
\nc\Pro{\text{\rm Prod}^{}_\alpha(F\op,\ab)}
\nc\Yt[1]{{\left.\Hom_\ct^{}\left(-,#1\right)\right|}_F^{}}
\nc\dl{\delta}
\nc\Proj[1]{#1\text{--\rm Proj}}
\nc\proj[1]{#1\text{--\rm proj}}
\nc\Flat[1]{#1\text{--\rm Flat}}
\nc\Inj[1]{#1\text{--\rm Inj}}
\nc\Ima{\mathrm{Im}}
\nc\Ker{\mathrm{Ker}}
\nc\ov{\overline}
\nc\wt{\widetilde}
\nc\wh{\widehat}
\nc\ph{\varphi}
\nc\tstr{{\it t}--structure}
\nc\tst[1]{\left({#1}^{\leq0},{#1}^{\geq0}\right)}
\nc\tstv[2]{\left({#1}_{#2}^{\leq0},{#1}_{#2}^{\geq0}\right)}
\nc\tsth[2]{{#1}_{#2}^{\heartsuit}}

\nc\spec[1]{{\text{\rm Spec}(#1)}}

\newcommand{\fc}{\mathfrak{C}}
\newcommand{\fl}{\mathfrak{L}}
\newcommand{\fs}{\mathfrak{S}}
\nc\EProd{\text{\rm EProd}}
\nc\ECoprod{\text{\rm ECoprod}}
\nc\Prod{\text{\rm Prod}}
\nc\ldimp{\text{\rm LDim}^{\prod}}
\nc\ldimc{\text{\rm LDim}^{\coprod}}
\nc\gen[2]{{\langle#1\rangle}^{}_{#2}}
\nc\genu[3]{{\langle#1\rangle}^{[#3]}_{#2}}
\nc\ogen[1]{\ov{\langle#1\rangle}}
\nc\ogenun[2]{\ov{\langle#1\rangle}_{#2}^{}}
\nc\ogenu[3]{\ov{\langle#1\rangle}^{[#3]}_{#2}}
\nc\ogenul[3]{\ov{\langle#1\rangle}^{(-\infty,#3]}_{#2}}
\nc\ogenuf[3]{\ov{\langle#1\rangle}^{[#3,\infty)}_{#2}}
\nc\genuf[3]{{\langle#1\rangle}^{[#3,\infty)}_{#2}}
\nc\genul[3]{{\langle#1\rangle}^{(-\infty,#3]}_{#2}}
\nc\dperf[1]{\D^{\mathrm{perf}}(#1)}
\nc\dcoh{\mathbf{D}^b_{\mathrm{coh}}}
\nc\dperfs[2]{\D_{#1}^{\mathrm{perf}}(#2)}
\nc\dcohs[1]{\mathbf{D}^b_{\mathrm{coh},#1}}

\newcommand{\Dqcs}[1]{{\mathbf D_{\text{\bf qc},#1}}}

\nc\dmcoh{\mathbf{D}^-_{\mathrm{coh}}}
\nc\dmcohs[1]{\mathbf{D}^-_{\mathrm{coh,#1}}}
\nc\dscoh{\mathbf{D}^{}_{\mathrm{coh}}}
\nc\RHHom{{\script{RH}}{\mathrm{om}}}
\nc\Coprod{\mathrm{Coprod}}
\nc\COprod{\mathrm{coprod}}
\nc\add{\mathrm{add}}
\nc\Add{\mathrm{Add}}
\nc\Smr{\mathrm{smd}}
\nc\id{\mathrm{id}}
\nc\LL{\mathbf{L}}
\nc\R{\mathbf{R}}
\nc\tsb{\ct^{\mathrm{sb}}}
\nc\Projdg[1]{\text{\rm Proj}^{}_{\mathbf{dg}}#1}

\renewcommand{\leq}{\leqslant}
\renewcommand{\geq}{\geqslant}

\nc\hoco{
\begin{picture}(40,10)
\put(20,0){\makebox(0,0)[b]{\text{\rm Hocolim}}}
\put(5,-2){\vector(1,0){30}}
\end{picture}\,\,}

\nc\holim{
\begin{picture}(40,10)
\put(20,0){\makebox(0,0)[b]{\text{\rm Holim}}}
\put(35,-2){\vector(-1,0){30}}
\end{picture}}

\begin{document}

\author{Amnon Neeman}\thanks{The research began at a workshop
  in Oberwolfach, and at the time was partly supported 
  by the Australian Research Council (grant number
  DP150102313). Major revisions were made recently and
  supported by the ERC Advanced Grant 101095900-TriCatApp.
  The author would like to thank
  all three institutions.}
\address{Dipartimento di Matematica ``F.\ Enriques''\\
        Universit{\`a} degli Studi di Milano\\
        Via Cesare Saldini 50\\
	20133 Milano\\
        ITALY}
\email{amnon.neeman@unimi.it}

\title{The categories $\ct^c$ and $\ct^b_c$ determine each other}

\begin{abstract}
  The main result is very general---it works in the abstract setting of weakly
  approximable triangulated categories. But it has the following concrete,
  immediate corollaries.
  \begin{enumerate}
    \item
      Suppose $X$ is a noetherian scheme. There is a recipe which, out of the category $\dperf X$, constructs $\dcoh(X)$ as a triangulated category.
    \item There is a recipe which, out of the category $\dcoh(X)$, constructs $\dperf X$ as a triangulated category.
    \item
      Let $R$ be any ring, possibly noncommutative. The recipe takes the triangulated category $\D^b(\proj R)$, that is the category whose
      objects are bounded complexes of finitely generated, projective left $R$-modules,
      and out of it constructs the triangulated category $\D^{-,b}(\proj R)$---that is the category of bounded-above cochain
      complexes of finitely generated projective left $R$-modules,
      with bounded cohomology.
    \item Now assume $R$ is left-coherent.
      Starting with $\D^b(\mod R)$ we construct $\D^b(\proj R)$.
    \item Out of the homotopy category of finite spectra we construct the homotopy category of spectra with finitely many nonzero stable homotopy groups, all of them finitely generated.
      \item 
      Out of the homotopy category of spectra with finitely many nonzero stable homotopy groups, all of them finitely generated, we construct the homotopy category of finite spectra.
  \end{enumerate}
More abstractly: given a category $\cs$ there is an old notion
of \emph{metrics} on it---this goes back to the 1970s,
see the articles by Lawvere~\cite{Lawvere73} and
Betti and Galuzzi~\cite{Betti-Galuzzi75}.
If $\cs$ is triangulated, we may require the metric to be
\emph{good,} in a sense that we will define in this article.
We may complete any essentially small triangulated category $\cs$ with
respect to any metric, obtaining a category $\fl(\cs)$ which isn't
usually triangulated. But inside $\fl(\cs)$ there is a subcategory
$\fs(\cs)$, of objects \emph{compactly supported} with respect to the metric.
And the first main theorem tells us that, as
long as the metric is good, the
category $\fs(\cs)$ is always triangulated.
The second main theorem gives a practical method
that can help in computing $\fs(\cs)$, as a triangulated category.

In the numbered
examples above, the metric on each of the triangulated categories
can be described intrinsically---it isn't added structure. There
are recipes
that start with essentially small
triangulated categories and, under some weak
hypotheses, cook up metrics.
\end{abstract}



\maketitle

\tableofcontents

\setcounter{section}{-1}

%\large
\section{Introduction}
\label{S0}

This article began with a question from Krause. In Oberwolfach,
in March 2018, Krause told me he had reread Rickard's old 
paper~\cite{Rickard89b} on derived
Morita equivalence, and didn't believe some of the proofs.
Specifically: Rickard~\cite[Theorem~6.4]{Rickard89b} 
asserts that, if $R$ and $S$ are coherent rings, then
\[
\D^b(\proj R)\cong \D^b(\proj S) \quad\Longleftrightarrow\quad
\D^b(\mod R)\cong \D^b(\mod S)
\]
Krause asked if I could find a counterexample.
More specifically: he was wondering
if the direction $\Longleftarrow$ is true.
Instead of a counterexample I discovered a proof, based on the ideas of
approximability---it's easily modified to be
symmetric enough to do both
directions. This article is about a vast generalization, but
for the sake of clarity it seems best to start with the simple Oberwolfach
argument. But first some notation.

\rmd{R0.1}
Suppose $\cs$ is a triangulated category, $G\in\cs$ is an
object and $A\leq B$ are integers. Then 
$\genu G{}{A,B}$ is the smallest full subcategory
of $\cs$, containing $\T^{-i}G$ for
$A\leq i\leq B$, and closed under direct summands and extensions.

In the above we allow $A$ and/or $B$ to be infinite. For example:
the case when $A=-\infty$ and $B=\infty$ gives a subcategory
${\langle G\rangle}_{}^{(-\infty,\infty)}$, containing all suspensions
of $G$. It is usually abbreviated $\gen G{}$, and is the smallest thick
subcategory containing $G$.

An object $G\in\cs$ is called a \emph{classical generator} if
$\cs=\gen G{}$.

Suppose $\cs$ is essentially small.
The category $\MMod\cs$ is the category of all additive functors
$H:\cs\op\la\ab$. The Yoneda functor $Y:\cs\la\MMod\cs$,
taking an object $A\in\cs$ to
$Y(A)=\Hom(-,A)\in\MMod\cs$, is a fully faithful embedding.

If $\cc$ is any pointed category and $\Cp\subset\cc$ is a subcategory,
then $\Cp^\perp$ is the
full subcategory of all objects $c\in\cc$ with $\Hom(\Cp,c)=0$.
And $^\perp\Cp$ is the full subcategory of all $c\in\cc$ with
$\Hom(c,\Cp)=0$. We will allow ourselves to take perpendiculars
in both $\cs$ and in $\MMod\cs$. To avoid confusion, when
we give a subcategory $\Cp\subset\cs$ we will write $\Cp^\perp$
for its perpendicular in $\cs$, and $Y(\Cp)^\perp$ for its perpendicular
in $\MMod\cs$.
\ermd

Now we come to something new:

\dfn{D0.3}
Let $\cs$ be a
triangulated category, and suppose we are given a
sequence of subcategories $\{\Cp_i\subset\cs,\,i\in\nn\}$ with
$\Tm\Cp_i\cup\Cp_i\cup\T\Cp_i\subset\Cp_{i+1}$. A sequence in $\cs$ of
the form $E_1\la E_2\la E_3\la\cdots$ is declared to be \emph{Cauchy
with respect to $\{\Cp_i\subset\cs,\,n\in\nn\}$} 
if, for every $i\in\nn$, there exists an integer $N>0$ such that
$\Hom(P,-)$ takes $E_n\la E_{n+1}$ to an isomorphism for all
$n\geq N$ and all $P\in\Cp_i$.
\edfn

\rmk{R0.5}
There is an obvious notion of equivalence---two sequences
of subcategories are declared equivalent if they yield the same
Cauchy sequences. For example: if $\cs$ has a classical generator
$G$ we can define $\Cp_i(G)=\genuf G{}{-i}$, and the resulting
Cauchy sequences don't depend on the choice of $G$.
The Cauchy sequences with respect to
this particular ``metric'' are
intrinsic, they depend only on $\cs$.
\ermk

\exm{E0.7}
If $\cs=\D^b(\proj R)$, with $R$ a ring, then the object $R$ is
a classical generator. Remark~\ref{R0.5} gives
an intrinsic notion of Cauchy sequences---to compute what
they are let us put $\Cp_i=\Cp_i(R)=\genuf R{}{-i}$ as
in Remark~\ref{R0.5}. The reader can check
that sequence $E_1\la E_2\la E_3\la\cdots$
is Cauchy precisely if, for every integer $i>0$, there exists
an integer $N>0$ such that $H^j(E_n)\la H^j(E_{n+1})$ is an isomorphism
whenever $n\geq N$ and $j\geq-i$.
\eexm

\dfn{D0.9}
Let $\cs$ be an essentially small triangulated category.
Let the notation be
as in Definition~\ref{D0.3}: that is $\{\Cp_i\subset\cs,\,i\in\nn\}$
is a sequence of
subcategories satisfying the hypotheses, with corresponding
Cauchy sequences.
We let $\fl(\cs)$
be the full subcategory of $\MMod\cs$
whose objects are the colimits
of Cauchy sequences in $\cs$. And we
declare
\[
\fs(\cs)\eq\fl(\cs)\cap\left[\bigcup_{i\in\nn}
  Y\big(\Cp_i^\perp\big)^\perp\right]
\]
\edfn

\exm{E0.11}
It is a small exercise to check that, in the case where $\cs=\D^b(\proj R)$
and the Cauchy sequences are as in Example~\ref{E0.7}, the category
$\fl(\cs)$ comes down to the image in $\MMod\cs$ of $\D^-(\proj R)$,
and if $R$ is coherent
the category $\fs(\cs)$ is nothing other than $\D^b(\mod R)$.
We have found an intrinsic way to construct $\D^b(\mod R)$ out
of $\D^b(\proj R)$.
\eexm

It is interesting to go in the other direction. Our problem
becomes to intrinsically define Cauchy sequences in $\D^b(\mod R)$. Note
that, except for special cases,
I have no idea when $\D^b(\mod R)$ has a classical
generator. When $R$ is commutative and noetherian
a great deal is known---the category $\D^b(\mod R)$
often has a classical generator, with the oldest theorem
being
Rouquier~\cite[Theorem~7.38]{Rouquier08}.
However: even for noetherian, commutative rings there
exist (pathological) counterexamples, where
$\D^b(\mod R)$ does not have a classical generator.
And here we're working in the generality of any, possibly noncommutative,
non-noetherian rings.

In the absence of a
classical generator the recipe of Remark~\ref{R0.5} isn't much use,
we need an alternative.

\dfn{D0.13}
Let $\cs$ be a triangulated category. We define a partial order on its
subcategories: we declare $\Cp\preceq\cq$ if there exists an
integer $n>0$ with $\T^n\Cp\subset\cq$. And $\Cp,\cq$ are declared
\emph{equivalent} if $\Cp\preceq\cq\preceq\Cp$.

For any object $G\in\cs$ consider the
subcategory $\big[\genul G{}0\big]^\perp$.
If, up to equivalence, there is a unique minimal
$\big[\genul G{}0\big]^\perp$ with respect to the partial order above,
we choose a member in this equivalence class and
call it $\cq(\cs)$.
\edfn

\exm{E0.15}
If $R$ is coherent and
$\cs=\D^b(\mod R)$ then, up to equivalence, the 
subcategories
$\big[\genul G{}0\big]^\perp$ have a unique minimal member.
After all: any object $G\in\D^b(\mod R)$ is contained in $\D^b(\mod R)^{\leq n}$
for some $n$, therefore $\big[\genul G{}0\big]^\perp$ contains
$\D^b(\mod R)^{\geq n+1}$. But if we take $G=R$ then
$\big[\genul G{}0\big]^\perp=\D^b(\mod R)^{\geq 1}$ and must therefore,
up to equivalence, be the
unique minimal subcategory.
\eexm

\dfn{D0.17}
Assume $\cs$ is an essentially small triangulated category, and assume
that a unique minimal $\cq(\cs)$ as in Definition~\ref{D0.13} exists.
Let the increasing sequence of subcategories $\{\Cp_i\subset\cs,\,i\in\nn\}$
be $\Cp_i=\T^i\cq(\cs)$.
\edfn

\exm{D0.19}
Now apply the construction of Definition~\ref{D0.9} to
$\cs=\big[\D^b(\mod R)\big]\op$ and to the sequence of categories
$\{\Cp_i\op\subset\cs,\,i\in\nn\}$, with $\Cp_i=\T^i\cq(\cs)$
as in Definition~\ref{D0.17}. A Cauchy
sequence turns out to be an inverse system
$\cdots\la E_3\la E_2\la E_1$ in $\D^b(\mod R)$, such that for every
$i>0$ there exists an $N>0$ with $H^j(E_{n+1})\la H^j(E_n)$ an
isomorphism whenever $n\geq N$ and $j\geq-i$. The category
$\cl(\cs)$ comes down to the image in $\MMod\cs$ of
the category $\big[\D^-(\mod R)\big]\op$,
and the category $\fs(\cs)$ is nothing other than
$\big[\D^b(\proj R)\big]\op$. We have found a recipe that goes back.
\eexm

We
have so far explained the simple
idea that led to this article.
Now it's time to elaborate on how we expand the ideas---it's time
to tell the reader what else she can expect to find in
the article, beyond the simple argument of the last couple of pages.
And the main issue comes down to:

\plm{P0.999333555666}
It is easy to reconstruct $\D^b(\mod R)$ out of $\D^b(\proj R)$
as a category---this much is a relatively short, straightforward
exercise, and we have outlined one solution above.
But Rickard's theorem says that, given rings
$R$ and $S$, the existence of a \emph{triangulated}
equivalence $\D^b(\proj R)\cong\D^b(\proj S)$ implies 
the existence of a \emph{triangulated}
equivalence $\D^b(\mod R)\cong\D^b(\mod S)$. In other
words the issue here is how to recover, out
of the \emph{triangulated} category $\D^b(\proj R)$, the
category $\D^b(\mod R)$ together with
\emph{its triangulated structure.}
\eplm

And this problem leads us
to the significant and surprising
results of the current article.
Before elaborating, let us address
one question I have been asked.

\rmk{R0.9993335556661}
In the last three decades there have
been stunning advances in the study of
enhancements of triangulated
categories---there has been beautiful
work by many authors, with the latest
and most powerful being the lovely
theory of (stable) infinity
categories. And it
is natural to wonder if this can be used
to help
with Problem~\ref{P0.999333555666}.

Let us explain this
as follows: we
could ask ourselves if, starting
with $\D^b(\proj R)$ \emph{together with
its structure as a stable $\infty$-category},
we could construct $\D^b(\mod R)$ \emph{with its
structure as a stable $\infty$-category.}

Answering this question in
any detail would take us much too
far afield. Let us confine ourselves
to referring to the fascinating
work of 
Denis-Charles
Cisinski and some
of his students, who have adapted
the methods of
the current article to the
stable $\infty$-category setting.
This leads
to a beautiful theory we will not
discuss here, beyond the comment
that the results of this article
are interesting also in their
stable $\infty$-category incarnation.

But there is another aspect,
to the
stable-infinity
enhancement story,
which we do want to elaborate
a tiny bit.
One of
the key points, of the current article,
is that the triangulated
structure we are about to create
is
\emph{strongly independent of the
stable $\infty$-category enhancement chosen.}
In sequels to this article, starting
with \cite{Neeman25},
we
will show that the work done here
can be developed into a powerful
new tool for proving unexpected
(strong) uniqueness-of-enhancement theorems.
To give but one application: a beautiful result of
Schwede~\cite{Schwede07}
showed that the homotopy category of
finite spectra has a unique
stable $\infty$-enhancement,
and using the techniques of this
article, together with Schwede's theorem,
allows one to produce a string of
other homotopy categories of spectra
with unique 
stable $\infty$-enhancements.
\ermk

Leaving the world of
enhancements,  let $\cs$ be any
essentially
small category. Modulo a slight change
of terminology, in the classical
literature one can find 
the notion of a
\emph{metric} on $\cs$, introduced long ago
in Lawvere~\cite{Lawvere73} and
Betti and Galuzzi~\cite{Betti-Galuzzi75}. The
terminological revision is that what we
call ``metric'' goes by the name ``norm''
in the old literature, and we reserve
the term ``norm'' for a very special class
of metrics, studied in other articles of
this series. 
The reader
can find much more detail in the survey
articles~\cite{Neeman19,Neeman22}.

So much for old stuff.
What is new
in this article is that we restrict
attention to \emph{triangulated categories,}
and, starting in Definition~\ref{D20.1}, we
focus only on \emph{good metrics.}
As in Definition~\ref{D0.9} we will define,
for any good metric on the triangulated category
$\cs$,
two full subcategories
\[
\fs(\cs)\subset\fl(\cs)\subset\MMod\cs\ .
\]
And the first theorem will be

\thm{T0.23}
For any essentially small triangulated category
$\cs$, and any good metric on $\cs$, the category $\fs(\cs)$ has a triangulated
structure which can be defined purely in terms of $\cs$ and the metric.
\ethm

\nin
We need hardly tell the reader how remarkable this is---there
are not many known recipes that start with a triangulated category
$\cs$, and out of it cook up another. The conventional wisdom is
that this can only be done in the presence of some enhancement.
But, in defiance of
conventional wisdom, in this article there is
no enhancement.

It becomes interesting to compute $\fs(\cs)$ in examples. For this it turns
out to be helpful to study the following situation.
\ntn{N0.25}
Let $\cs$ be an essentially small triangulated category with a metric.
Suppose we are given a
fully faithful triangulated functor $F:\cs\la\ct$; we consider
also the functor $\cy:\ct\la\MMod\cs$, which takes an object $A\in\ct$
to the functor $\Hom\big(F(-),A\big)$. The functor $F$ is called a
\emph{good extension with respect to the metric} if $\ct$ has
enough coproducts to define the homotopy colimit
$\hoco F(E_*)$ for any Cauchy sequence
$E_*$ in $\cs$,
and for every Cauchy sequence $E_*$ in $\cs$ the natural map
$\colim\,Y(E_*)\la \cy\big(\hoco F(E_*)\big)$ is an isomorphism.
\entn

For any good extension $F:\cs\la\ct$ we proceed to define
the full subcategory  $\fl'(\cs)\subset\ct$ to have for objects
all the homotopy colimits of Cauchy sequences,
and inside $\fl'(\cs)$
we define a full subcategory $\fl'(\cs)\cap\cy^{-1}\big(\fc(\cs)\big)$---we
ask the reader for patience, the definition will come
in the body of the
paper\footnote{More precisely: for the category $\fl'(\cs)$ see Definition~\ref{D21.7}, for the functor $\cy$ see Notation~\ref{N21.-100}, and the category $\fs(\cs)$ is introduced in Definition~\ref{D20.11}(2).}.
The next result is

\thm{T0.27}
The category $\fl'(\cs)\cap\cy^{-1}\big(\fc(\cs)\big)$ is a
triangulated subcategory of $\ct$, and the functor
$\cy:\ct\la\MMod\cs$ restricts to a
triangulated equivalence
\[\xymatrix@C+30pt{
\cy\,\,:\,\,\fl'(\cs)\cap\cy^{-1}\big(\fc(\cs)\big)\ar[r] &\fs(\cs)
}\]
\ethm

\nin
In the presence of a good extension $F:\cs\la\ct$ this allows us to
compute $\fs(\cs)$, up to triangulated equivalence, as the
triangulated subcategory $\fl'(\cs)\cap\cy^{-1}\big(\fc(\cs)\big)$ of $\ct$.

Next suppose $\ct$ is a triangulated category with coproducts, and
assume it has a compact generator $H$ with $\Hom(H,\T^iH)=0$ for $i\gg0$.
In this case the theory introduced in \cite{Neeman17,Neeman24}
kicks in:
there is a \emph{preferred equivalence class} of {\it t}--structures
on $\ct$, and it is possible to define, \emph{intrinsically,}
thick subcategories $\ct^b_c\subset\ct^-_c$. In terms of the
preferred {\it t}--structures it is possible to endow $\ct^c$ and
$\big[\ct^b_c\big]\op$
with metrics, we will define them in Example~\ref{E20.3}.
It turns out that the embedding $\ct^c\la\ct$ is always a good
extension, while the embedding $\big[\ct^b_c\big]\op\la\ct\op$
is a good extension provided $\ct$ is weakly approximable as in \cite[Definition~0.25]{Neeman24}.
And the next result says that, in full, gorgeous generality

\pro{P0.29}
Assume that category $\ct$ is weakly approximable.
Then for the metrics above we have
\be
\item
$\fs(\ct^c)=\ct^b_c$.
\item
If $\ct$ is coherent and weakly approximable then
$\fs\big(\big[\ct^b_c\big]\op\big)=\big[\ct^c\big]\op$.
\ee
\epro
\nin
The notion of a \emph{coherent} triangulated category,
in Proposition~\ref{P0.29}(ii), is new. It will
be defined
in Section~\ref{S29}. It is a hypothesis that guarantees there
are enough nonzero objects in $\ct^b_c$, after all there is no \emph{a priori}
reason to expect any.

From the perspective of the Oberwolfach discussion this is still
unsatisfactory: in the special case where $\ct=\D(\Mod R)$, the derived
category of a coherent ring, the recipe tells us how to
pass from
\[\begin{array}{ccc}
\ct^c=\D^b(\proj R),\text{ together with its metric} &\Longrightarrow & \D^b(\mod R)=\ct^b_c\\ 
\big[\ct^b_c\big]\op=\D^b(\mod R)\op,\text{ together with its metric} &\Longrightarrow & \D^b(\mod R)\op=\big[\ct^b_c\big]\op\ .
\end{array}\]
But the metrics are defined in terms of the
preferred equivalence class of {\it t}--structures on $\ct$. As presented,
the metrics depend
on the embedding into $\ct$.

Hence it becomes interesting to see when we can construct the
metrics intrinsically, without reference to $\ct$. It turns out we can
always do this. More precisely: the recipe of Definitions~\ref{D0.13}
and \ref{D0.17} works in general, to give the metric on $\big[\ct^b_c\big]\op$
in Proposition~\ref{P0.29}(ii).
There is also an intrinsic description of the metric on $\ct^c$
used in Proposition~\ref{P0.29}(i), we will see it
in Definition~\ref{D22.103}(i) and Remark~\ref{R22.105.5}, but it
isn't the recipe given in Remark~\ref{R0.5}. For the metric of
Remark~\ref{R0.5} to agree with the metric
in Proposition~\ref{P0.29}(i) we will
need to assume $\ct$ weakly approximable,
see Proposition~\ref{P22.109}.

\rmk{R0.99937}
Let $X$ be a noetherian scheme, and let $Z\subset X$
be a Zariski-closed subset. Let $\Dqcs Z(X)$ denote
the derived category, whose objects $F\in\Dqcs Z(X)$
are complexes
of $\co_X^{}$--modules satisfying two hypotheses:
\be
\item
The cohomology sheaves $\ch^i(F)$ are all quasicoherent.
\item
The restriction of $F$ to the open set $X-Z$ is acyclic.
\ee
From \cite[Theorem~3.2(iv)]{Neeman22A} we learn that
the category $\Dqcs Z(X)$ is weakly approximable,
and it is easy to see that it is coherent.

Proposition~\ref{P0.29} therefore applies, and tells
us that, with the notation meant to be self-explanatory,
\[
\fs\left(\dperfs ZX\right)=\dcohs Z(X)
\qquad\text{and}\qquad
\fs\left(\dcohs Z(X)\op\right)={\dperfs ZX}\op\ .
\]
For the sake of simplicity let us specialize to
the absolute case, where $Z=X$. The above specializes
to
\[
\fs\left(\dperf X\right)=\dcoh (X)
\qquad\text{and}\qquad
\fs\left(\dcoh(X)\op\right)={\dperf X}\op\ .
\]
In other words: this article gives an
explicit recipe which, for any noetherian scheme $X$,
passes back and forth between the triangulated categories
$\dperf X$ and $\dcoh(X)$.

To underline just how surprising this is let me
remind the reader that, for at least two
decades now, birational geometers have been
using derived categories as invariants of
noetherian schemes. The philosophy is that, in the
minimal model program, flips conjecturally make
the derived category ``smaller'' while flops
leave the derived category almost
unchanged\footnote{The ``almost'' in this context means
that, conjecturally, there is only a finite
set of flop ``partners'' of any projective variety.}.

Now recall: if the dimension of $X$ is $\geq3$
then the minimal models will inevitably be
singular---let's not enter into a discussion of
which terminal singularities we are going to
permit, the point we want to make is that
the categories $\dcoh(X)$ and $\dperf X$ are
\emph{not going to agree.} And birational geometers
have been known to argue over which is the right
invariant.

From this article we learn that this particular
argument is a waste of time: the categories
$\dcoh(X)$ and $\dperf X$ are interchangeable,
each of them contains all the information about the
other.
\ermk

\rmk{R0.2121029}
Returning to Oberwolfach in 2018:
Krause asked me the question that
led to this article, but he also
went on to work on it independently.
And for the 
direction of passing from $\D^b(\proj R)$ to $\D^b(\mod R)$,
Krause~\cite{Krause18}
developed a different argument.
Without going into detail, let
us confine ourselves to saying that
the results presented in this
article are far better and
more general.
\ermk

\rmk{R0.99}
In Remark~\ref{R0.2121029} we mentioned that the
development in Krause~\cite{Krause18} is
different.

In Krause's article there is
no mention of metrics, he speaks only
of ``Cauchy sequences''. And he only considers
one class of Cauchy sequences. Hence one
of the innovations of the current
article is to view the Cauchy sequences
as coming from metrics, and permit ourselves
to study a variety of metrics.

Stated this way,
the general theory developed
in this article applies to Krause's
Cauchy sequences---they do come from a choice of
metric. Krause's metric
happens to be different from mine.
It is of course possible to apply
to Krause's metric
the methods developed in this
article.
I have 
only fully worked out, in detail,
what this gives
in the special case where $\cs=\D^b(\proj R)$, with
$R$ a noetherian ring. In
this case the triangulated category $\fs(\cs)$ turns out
to be $\D^b_{\mod R}(\Inj R)$, the category of bounded complexes
of injective $R$--modules whose cohomology modules are finite.
If $R$ has a dualizing complex this category is
equivalent to $\D^b(\proj R)$. For more detail see
Examples~\ref{E21.-5} and \ref{E21.10},
as well as
Remark~\ref{R22.856}.
\ermk

\rmk{R0.212121}
There is a philosophical point here that deserves underlining.

The best way to look at the current manuscript is in
the context of approximable triangulated categories,
a notion introduced in \cite{Neeman24}, which
we already mentioned in the leadup
to Proposition~\ref{P0.29}. And at the core
of the theory is the idea of borrowing techniques from
real analysis, and adapting and modifying them to
work in the world of triangulated categories. And given this,
we should learn from the experience of our analyst friends.
And their experience teaches us that making
a wise choice of metrics is crucial for solving the
problem at hand.

After all normally, when we set about solving a
partial differential equation, we begin with a weak
solution---usually the weak solution is easy to find, and
it lies in some humongous space of distributions.
The hard work is in proving regularity, that is in showing
(for example) that the solution is $C^\infty$.

And invariably the proof involves estimating some
Sobolev norms, and combining these estimates using
a string of Sobolev inequalities.

For any given problem there will be a choice of metrics
that gives elegant proofs of strong results, and then
there will be many poor choices of metrics. In this article
we will attempt to show that the metric we choose is the
natural one for the problem at hand. It delivers powerful
results, does so elegantly and in great generality.
\ermk

Finally we should say something about the structure of the article.
The first two sections work with a triangulated category $\cs$ and
its metric---there is no mention of good extensions, the sections
are devoted to the proof of Theorem~\ref{T0.23} and are
self-contained. 
Section~\ref{S21}
is where we prove Theorem~\ref{T0.27}---if the reader ignores the
examples, Section~\ref{S21} is also self-contained. But the later sections, which
work out the general theory in the examples $\ct^c$ and $\big[\ct^b_c\big]\op$,
assume familiarity with approximability.

\medskip

\nin
{\bf Acknowledgements.}\ \ The author would like to thank Henning Krause for
suggesting the problem, and the Mathematisches Forschungsinstitut Oberwolfach
for its hospitality and congenial working environment during the week
when the work started, back in March 2018.
Further thanks go to
Alberto Canonaco, Patrick Lank and Paolo Stellari
for corrections and improvements to earlier versions.

\section{The basic definitions}
\label{S20}

\rmd{R-500.1}
Let $\cs$ be a triangulated category and let
$\ca,\cc$ be subcategories. As in
\cite[1.3.9]{BeiBerDel82}
we define the full subcategory $\ca*\cc\subset\cs$ to have for objects
those
$b\in\cs$ for which there exists, in $\cs$,
a triangle $a\la b\la c\la$ with $a\in\ca$
and $c\in\cc$.
\ermd

\dfn{D20.1}
Let $\cs$ be a triangulated category. A \emph{good metric} on $\cs$ will be a
sequence
of additive subcategories $\{\cm_i\subset\cs\mid i\in\nn\}$,
containing $0$ and such that
\be
\item
$\cm_i$ contains  
the union
$\Tm\cm_{i+1}\cup\cm_{i+1}\cup\T\cm_{i+1}$
for every $i$.
\item
$\cm_i*\cm_i=\cm_i$.  
\ee
A good metric $\{\cm_i\}$ is declared to be \emph{finer than the good metric}
$\{\cn_i\}$
if, for every integer $i>0$, there exists an integer $j>0$ with
$\cm_j\subset\cn_i$; we denote this partial order by
$\{\cm_i\}\preceq\{\cn_i\}$.
The good metrics $\{\cm_i\}$, $\{\cn_i\}$ are \emph{equivalent} if
$\{\cm_i\}\preceq\{\cn_i\}\preceq\{\cm_i\}$.
\edfn

\exm{E20.2}
The dumb example is to let $\cm_i=\cs$ for every $i$.
\eexm

\rmd{R20.2.5}
For the next example we remind the reader of some constructions
from \cite{Neeman24}. Let $\ct$ be a triangulated category.
In \cite[Definition~0.14]{Neeman24} we declared two
{\it t}--structures $\big(\ct_1^{\leq0},\ct_2^{\geq0}\big)$ and
$\big(\ct_2^{\leq0},\ct_2^{\geq0}\big)$ to be \emph{equivalent} if
there exists an integer $A>0$ with
$\ct_1^{\leq-A}\subset\ct_2^{\leq0}\subset\ct_1^{\leq A}$.
Now assume $\ct$ has coproducts and a compact generator $G$.
Then \cite[Definition~0.18]{Neeman24} defines a \emph{preferred
  equivalence class} of {\it t}--structures, namely the one
containing the \tstr\ generated by $G$ in the sense
of Alonso, Jerem\'{\i}as and Souto~\cite{Alonso-Jeremias-Souto03}.
And \cite[Definition~0.20]{Neeman24} allows one to construct
two subcategories $\ct^b_c\subset\ct^-_c\subset\ct$. If the
compact generator $G\in\ct$ is such that $\Hom(G,\T^iG)=0$ for $i\gg0$,
then $\ct^b_c\subset\ct^-_c$ are both thick subcategories of $\ct$,
see \cite[Proposition~2.10]{Neeman24}.
\ermd

\exm{E20.3}
Suppose $\ct$ is a triangulated category with coproducts, and assume
$\ct$ has a compact generator $G$ with $\Hom(G,\T^iG)=0$ for $i\gg0$.
With the notation as in Reminder~\ref{R20.2.5}
let $\big(\ct^{\leq0},\ct^{\geq0}\big)$ be a \tstr\ in the
preferred equivalence
class.
Out of this data we can construct two examples of $\cs$'s with good metrics:
\be
\item
Let $\cs$ be the subcategory $\ct^c\subset\ct$, and
put $\cm_i=\ct^c\cap\ct^{\leq-i}$.
\item
Let $\cs$ be the subcategory $\big[\ct^b_c\big]\op$, and put
$\cm\op_i=\ct_c^b\cap\ct^{\leq-i}$.
\ee
It's obvious that equivalent {\it t}--structures define equivalent good metrics.
Thus up to equivalence we have a canonical good metric on $\ct^c$
and a canonical good metric on $\big[\ct^b_c\big]\op$. But the definition
depends on the embedding into $\ct$, which is the category
with the \tstr.
\eexm

\dfn{D20.5}
Let $\cs$ be a triangulated category with a good metric $\{\cm_i\}$.
A \emph{Cauchy sequence} in $\cs$ is a sequence of
$E_1\la E_2\la E_3\la\cdots$ so that, for every
integer $i>0$,
there exists
an integer $M>0$ such that,  in any triangle 
$E_m\la E_{m'}\la D_{m,m'}$ with $M\leq m<m'$, the object $D_{m,m'}$ lies
in $\cm_i$.  
\edfn

\rmk{R20.6}
In Definition~\ref{D20.1}(i) we postulated
that $\cm_i$ contains the union
$\Tm\cm_{i+1}\cup\cm_{i+1}\cup\T\cm_{i+1}$. 
An easy induction gives that
\[
\cm_i\quad\supset\quad\bigcup_{j=-n}^n\T^j\cm_{i+n}\ .
\]
It immediately follows that, given any Cauchy sequence
$E_1\la E_2\la E_3\la\cdots$ and any pair of
positive integers $\{i,n\}$,
there exists
an integer $M>0$ such that,  in any triangle 
$E_m\la E_{m'}\la D_{m,m'}$ with $M\leq m<m'$, the object $\T^jD_{m,m'}$ lies
in $\cm_i$ for all $-n\leq j\leq n$.  
\ermk

\rmk{R20.7}
Note that the Cauchy sequences depend only on the equivalence
class of the good metric. 
\ermk

The following observation is useful for constructing Cauchy sequences.

\lem{L20.8}
Suppose we are given in $\cs$ a sequence 
$E_1\la E_2\la E_3\la\cdots$ so that, for any integer $i>0$, 
there exists
an integer $M>0$ such that,  in any triangle 
$E_m\la E_{m+1}\la D_{m}$ with $M\leq m$, the object $D_{m}$ lies
in $\cm_i$.

Then
for all integers
$M\leq m< m'$, the triangles $E_m\la E_{m'}\la D_{m,m'}$ have
$D_{m,m'}\in\cm_i$. In particular: the sequence $E_*$ is Cauchy.
\elem

\prf
We prove the assertion by induction on $m'-m$, the case $m'-m=1$
being the hypothesis. Suppose the result is true for $m'-m\leq k$,
and construct an octahedron
on the composable maps $E_m\la E_{m+k}\la E_{m+k+1}$; we
obtain a triangle
$D_{m,m+k}\la D_{m,m+k+1}\la D_{m+k}$, and
induction gives that $D_{m,m+k},D_{m+k}\in\cm_i$.
It follows that $D_{m,m+k+1}\in\cm_i*\cm_i=\cm_i$.
\eprf

\rmk{R20.9}
Given an essentially small
triangulated category $\cs$ it is customary to consider the Yoneda
functor on it, we wish to explore the functor $Y:\cs\la\MMod\cs$.
To recall the notation: $\MMod\cs$ is the
category of additive functors $\cs\op\la\Mod\zz$, and the functor $Y$ takes
the object $A\in\cs$ to the additive functor $Y(A)=\Hom(-,A)$.
\ermk

\dfn{D20.11}
Suppose $\cs$ is an essentially small
triangulated category with a
good metric $\{\cm_i\}$.
We define three full subcategories $\fl(\cs)$, $\fc(\cs)$
and $\fs(\cs)$ of the category $\MMod\cs$ as follows:
\be
\item
The objects of $\fl(\cs)$ are
the functors $\cs\op\la\Mod\zz$ which can be
expressed as $\colim\, Y(E_i)$, where $E_*$ is a 
Cauchy sequence in $\cs$.
\item
The objects of $\fc(\cs)$ are the
compactly supported functors with respect to the
good metric.
Concretely, they are given by the formula
\[
\fc(\cs)\eq\left\{
A\in\MMod\cs\left|
\begin{array}{c}\text{There exists an integer }i>0\\
 \text{with }\Hom\big(Y(\cm_i)\,,\,A\big)=0
\end{array}\right.\right\}.
\]
\item
Finally, $\fs(\cs)$ is defined by
$\fs(\cs)=\fl(\cs)\cap\fc(\cs)$.
\ee
\edfn

\rmk{R20.13}
First of all: it's obvious that the categories $\fl(\cs)$, $\fc(\cs)$
and $\fs(\cs)$ depend only on the equivalence class
of the good metric.

Next note that the Yoneda functor $Y:\cs\la\MMod\cs$ factors through the
subcategory $\fl(\cs)$, after all the constant sequence
$E\stackrel\id\la E\stackrel\id\la E\stackrel\id\la\cdots$ is Cauchy
for any good metric, and the colimit in $\MMod\cs$ of the image of
this sequence under Yoneda is $Y(E)$.

Finally
observe that all the objects of $\fl(\cs)$ are homological functors
$\cs\op\la\Mod\zz$. After all they are filtered colimits of the homological
functors $Y(E_i)$.
\ermk

\exm{E20.1098}
In  the special case where the good metric is the dumb one in Example~\ref{E20.2},
that is $\cm_i=\cs$ for every $i$, every sequence is Cauchy and
$\fl(\cs)$ is the Ind-completion of $\cs$. The category $\fc(\cs)$ and
$\fs(\cs)$ are both equal to $\{0\}$. The theory doesn't produce much.
\eexm

\section{The category $\fs(\cs)$ is triangulated}
\label{S28}

\ntn{N28.1}
Throughout this section we will fix the triangulated category $\cs$ together
with its good metric $\{\cm_i\}$. The only categories we will study in the section
are full
subcategories of $\MMod\cs$: the subcategories $\fl(\cs)$, $\fc(\cs)$ and
$\fs(\cs)$ of Definition~\ref{D20.11}, as well as the subcategory
$\cs\subset\fl(\cs)$.
Note that we view $\cs$ as embedded in $\cl(\cs)\subset\MMod\cs$
through the fully faithful
functor $Y$. And most of the time we will freely confuse $\cs$ with
its image under $Y:\cs\la\MMod\cs$.

It's only in the statements of
results that we plan to appeal
to in later sections---not their
proofs---that
we try to be careful with the notation. The reason is that in later
sections we will allow ourselves to embed $\cs$ into other categories
$\ct$, and confusion could arise.
\entn

\dis{D28.3}
We define an invertible automorphism
$\T:\MMod\cs\la\MMod\cs$ by the rule
\be
\item
If $A$ is an object of $\MMod\cs$, meaning a functor $A:\cs\op\la\Mod\zz$,
and $s$ is an object of $\cs$, then
$[\T A](s)=A(\Tm s)$.
\setcounter{enumiv}{\value{enumi}}
\ee
The Yoneda isomorphism $A(s)\cong \Hom_{\MMod\cs}^{}\big(Y(s),A\big)$ permits
us, in our sloppy conventions established in
Notation~\ref{N28.1}, to rewrite (i) as $\Hom(s,\T A)=\Hom(\Tm s,A)$.
The homological functors $A:\cs\op\la\MMod \zz$ are precisely the objects
$A\in\MMod\cs$ such that $\Hom_{\MMod\cs}^{}(-,A)$ restricts to a homological
functor on $\cs\op$. In particular: in Remark~\ref{R20.13} we noted
\be
\setcounter{enumi}{\value{enumiv}}
\item
If $A$ belongs to $\cl(\cs)$, then the restriction to $\cs\op$ of the
functor $\Hom_{\MMod\cs}^{}(-,A)$ is homological.
\setcounter{enumiv}{\value{enumi}}
\ee
Finally, a sequence $A\la B\la C$ in $\MMod\cs$ is exact
if it is exact when evaluated at each $s\in\cs$. Our notation translates
this into
\be
\setcounter{enumi}{\value{enumiv}}
\item
The sequence $A\la B\la C$ in $\MMod\cs$ is exact
if and only if, for every object $s\in\cs$, the functor
$\Hom(s,-)$ takes it to an exact sequence.
\setcounter{enumiv}{\value{enumi}}
\ee
\edis

\obs{O28.-1}
Since the formula will be cited in future sections our notation is careful,
we write
\[
\fc(\cs)=
\bigcup_{i\in\nn}^{}\,\,\big[Y(\cm_i)\big]^\perp\ .
\]
Reverting to the sloppy conventions of
Notation~\ref{N28.1}
and Discussion~\ref{D28.3} for the explanation:
for $A\in\MMod\cs$ the condition
$\Hom_{\MMod\cs}^{}(\cm_i,A)=0$ rewrites
as $A\in\big[\cm_i\big]^\perp$, and the
displayed formula
above just codifies the quantifier on $i$ in
Definition~\ref{D20.11}(ii).

The above makes it clear that
\be
\item
$\fc(\cs)$ is closed in $\MMod\cs$ under direct summands.
\item
The fact that $\cm_i$ contains $\Tm\cm_{i+1}\cup\cm_{i+1}\cup\T\cm_{i+1}$
implies that  $\T\fc(\cs)=\fc(\cs)$.
\item
Given an exact sequence
$A\la B\la C$ in $\MMod\cs$, we have
\[
A,C\in\fc(\cs)\Longrightarrow B\in\fc(\cs)\ .
\]
\ee
Perhaps we should explain (iii).
Since $\cm_i,\cm_j$ both contain $\cm_{i+j}$ we have
\[
A\in\cm_i^\perp\quad\text{and}\quad C\in\cm_j^\perp\quad
\Longrightarrow\quad B\in\cm_{i+j}^\perp\ ,
\]
and hence
\[
A\in\bigcup_i\cm_i^\perp\quad\text{and}\quad C\in\bigcup_i\cm_i^\perp
\quad\Longrightarrow\quad B\in\bigcup_i\cm_{i}^\perp\ .
\]
\eobs

\dfn{D28.7}
Let $\cs\subset\MMod\cs$ be as in Notation~\ref{N28.1}. A
\emph{pre-triangle} is a diagram 
$A\stackrel{f}\la B\stackrel{g}\la C\stackrel{h}\la \T A$, in
the abelian category $\MMod\cs$, such that
$A\stackrel{f}\la B\stackrel{g}\la C\stackrel{h}\la \T A\stackrel{\T f}\la \T B$ is an exact sequence.
\edfn

\rmk{R28.109}
From Observations~\ref{O28.-1}~(ii) and (iii) we learn that, if
$A\stackrel{f}\la B\stackrel{g}\la C\stackrel{h}\la \T A$ is a
pre-triangle in $\MMod\cs$ and if two of $A,B,C$ lie in $\fc(\cs)$,
then so does the third.
\ermk

\rmk{R28.905}
Given any distinguished triangle
$a\stackrel{f}\la b\stackrel{g}\la c\stackrel{h}\la \T a$
in the category $\cs$, its
image under Yoneda
$Y(a)\stackrel{Y(f)}\la Y(b)\stackrel{Y(g)}\la Y(c)\stackrel{Y(h)}\la \T Y(a)$
is a pre-triangle in $\MMod\cs$. This comes from
combining the fact that $Y:\cs\la\MMod\cs$ is fully faithful,
with observation
that, for $s\in\cs$, the functor $\Hom(s,-)$
is a homological functor $\cs\la\ab$.

But now any filtered colimit of
pre-triangles is obviously a pre-triangle, and
(for us) an
important special case will be:
\ermk

\dfn{D28.907}
A \emph{strong triangle} in $\fl(\cs)$ is a sequence
$A\stackrel{f}\la B\stackrel{g}\la C\stackrel{h}\la \T A$
which is isomorphic to the colimit
in $\MMod\cs$ of the image under Yoneda of a Cauchy
sequence of distinguished triangles.
\edfn

Having made the definition, 
it becomes interesting to construct examples.

\lem{L28.19}
Let $f:A\la B$ be a morphism in $\fl(\cs)$. 
First of all, it may be completed to a
strong triangle 
$A\stackrel{f}\la B\stackrel{g}\la C\stackrel{h}\la \T A$ in the
category
$\fl(\cs)$.

But the following, more
precise assertions are also true:
\be
\item
Suppose we are given in $\cs$ any Cauchy sequences $a_*$ and $b'_*$
such that $A=\colim\, Y(a_*)$
and
$B=\colim\,Y(b'_*)$.
We may choose a subsequence $b_*$ of $b'_*$,
as well as a map of sequences $f_*:a_*\la b_*$, with $f=\colim\,Y(f_*)$.
Moreover: if we are given in advance a
morphism $f_1:a_1\la b'_1$ such that the square below commutes
\[\xymatrix@C+30pt{
Y(a_1^{})\ar[r]^-{Y(f'_1)}\ar[d]& Y(b'_1)\ar[d]\\
A\ar[r]^-f & B
}\]
then we may choose the subsequence $b_*\subset b'_*$ 
so that $b_1^{}=b'_1$ and $f_1=f'_1$.
\item
Given a Cauchy sequence of maps $f_*:a_*\la b_*$ in $\cs$,
we may complete it in $\cs$ to a Cauchy sequence of  
triangles
$a_*\stackrel{f_*}\la b_*\stackrel{g_*}\la c_*\stackrel{h_*}\la \T a_*$.
If the triangle 
$a_1^{}\stackrel{f_1}\la b_1^{}\stackrel{g_1^{}}\la c^{}_1\stackrel{h^{}_1}\la \T a_1^{}$
is specified in advance, we may stick with
it as the first triangle in
our Cauchy sequence
$a_*\stackrel{f_*}\la b_*\stackrel{g_*}\la c_*\stackrel{h_*}\la \T a_*$.
Applying the Yoneda functor $Y$
to the Cauchy sequence
and taking colimits produces
a strong triangle 
$A\stackrel{f}\la B\stackrel{g}\la C\stackrel{h}\la \T A$ in the
category
 $\fl(\cs)$.
\item
Finally, with $f:A\la B$ as in the statement of the
the Lemma and with $a_*$ and $b'_*$ as in (i), then
for any choice of 
map of sequences $f_*:a_*\la b_*$ as in (i),
the completion to
a Cauchy sequence of  
triangles
$a_*\stackrel{f_*}\la b_*\stackrel{g_*}\la c_*\stackrel{h_*}\la \T a_*$
in (ii) may be done in such a way that, for every pair
of integers $\ell,N>0$ 
for which the original 
Cauchy sequences $a_*$ and $b'_*$ of (i)
satisfy the hypothesis 
that if
$m<m'$ are any
integers, both $>N$, then in the triangles
$a_m\la a_{m'}\la \ov d_{m,m'}$ and
$b'_m\la b'_{m'}\la d_{m,m'}$ we have that
\[ \ov d_{m,m'},\T\ov d_{m,m'},
\Tm d_{m,m'}, d_{m,m'}\qquad\text{ all lie in }\qquad\cm_\ell\ ,
\]
then the triangles
$c_m\la c_{m'}\la \wh d_{m,m'}$ are such that $\Tm\wh d_{m,m'},\wh d_{m,m'}$
will also lie in $\cm_\ell$.
\ee
\elem

\prf
Because $A,B$ belong to $\fl(\cs)$ we can find Cauchy sequences
$a_*,b'_*\subset\cs$ with
$A\cong\colim\, a_*$ and $B\cong\colim\, b'_*$.  
If the sequences are given, as
in (i), we work with
those. But in any case we start with Cauchy sequences $a_*,b'_*\subset\cs$ with
$A\cong\colim\, a_*$ and $B\cong\colim\, b'_*$.  
We are given in $\cl(\cs)$ the composite $a_1^{}\la A\stackrel{f}\la B$,
which is an element in
$\colim\,\Hom(a_1^{},b'_\ell)$. We may choose a preimage in some
$\Hom(a_1^{},b'_\ell)$, constructing a commutative square
\[\xymatrix@C+20pt{
a_1^{}\ar[r]^-{f_1}\ar[d]& b'_{\ell_1^{}}\ar[d]\\
A\ar[r]^-f & B
}\]
If we are given $f_1$ we begin with it.

And then we continue inductively. If we have carried out the construction
as far as the integer $i$, then we have a commutative diagram
\[\xymatrix@C+20pt{
a_i\ar[r]^-{f_i}\ar[d]& b'_{\ell_i^{}}\ar[dd]\\
a_{i+1}\ar[d]     &     \\    
A\ar[r]^-f & B
}\]
which we may, by choosing $\ell_{i+1}^{}$ large enough, complete to a
commutative diagram
\[\xymatrix@C+20pt{
a_i\ar[r]^-{f_i}\ar[d]& b'_{\ell_i^{}}\ar[d]\\
a_{i+1}\ar[d]\ar[r]^-{f_{i+1}} &  b'_{\ell_{i+1}^{}} \ar[d] \\    
A\ar[r]^-f & B
}\]
constructing the subsequence $b_*$ of $b'_*$ and the map $f_*:a_*\la b_*$
This completes the proof of (i).

Next complete each $f_i:a_i\la b_i$ to a triangle
$a_i\stackrel{f_i}\la b_i\stackrel{g_i}\la c_i\stackrel{h_i}\la \T a_i$;
if the triangle
$a_1^{}\stackrel{f_1^{}}\la b_1^{}\stackrel{g_1^{}}\la c_1^{}\stackrel{h_1^{}}\la \T a_1^{}$
is already given we work with it.
In the paragraph above we constructed commutative squares
\[\xymatrix@C+20pt{
a_i\ar[r]^-{f_i}\ar[d]^{\wt\alpha_i} & b_i\ar[d]^{\wt\beta_i} \\
a_{i+1}\ar[r]^-{f_{i+1}} & b_{i+1}
}\]
which we complete to a $3\times3$ diagram of triangles
\[\xymatrix@C+20pt{
a_i\ar[r]^-{f_i}\ar[d]^{\wt\alpha_i} & b_i\ar[d]^{\wt\beta_i}\ar[r]^-{g_i} &
c_i\ar[r]^-{h_i}\ar[d]^{\wt\gamma_i} & \T a_i\ar[d]^{\T\wt\alpha_i} \\
a_{i+1}\ar[r]^-{f_{i+1}}\ar[d] & b_{i+1} \ar[r]^-{g_{i+1}}\ar[d] &
c_{i+1}\ar[r]^-{h_{i+1}}\ar[d] & \T a_{i+1}\ar[d]\\
\ov d_{i}\ar[r] & d_{i} \ar[r] &
\wh d_{i}\ar[r] & \T \ov d_{i} 
}\]
We note that, so far
\be
\item
We have extended the sequence of maps $f_*:a_*\la b_*$ to a
sequence of triangles
$a_*\stackrel{f_*}\la b_*\stackrel{g_*}\la c_*\stackrel{h_*}\la \T a_*$.
\setcounter{enumiv}{\value{enumi}}
\ee
Now let's step back
to the situation in (i): we are given
the sequences $a_*$, $b'_*$ are Cauchy
sequences with $A=\colim\,Y(a_*)$
and $B=\colim\, Y(b'_*)$. Of course:
if we are already in the situation as in (ii),
let $b'_*=b_*$.
Anyway: because $a_*$ and $b'_*$ are Cauchy,
for every integer $\ell>0$ we may choose
an integer $N>0$ such that, if $m<m'$ are integers
both $>N$, then 
in the triangles
$a_m\la a_{m'}\la \ov d_{m,m'}$ and
$b'_m\la b'_{m'}\la d_{m,m'}$ we have that
$\ov d_{m,m'},\T\ov d_{m,m'},\Tm d_{m,m'}, d_{m,m'}$
all lie in $\cm_\ell$.
Because $b_*$ is chosen to be a subsequence
of $b'_*$ we know, for all $m>0$, that $b_m=b'_n$ for some
integer $n\geq m$. Hence, with the same
integers $\ell,N$, we deduce that if $i$ is
any integer $>N$ then the  
triangles $a_i\la a^{}_{i+1}\la  d_i$
and $b_i\la b_{i+1}\la \ov d_i$,
of the $3\times3$ diagram of triangles
above,
satisfy $\Tm d_i,d_i,\ov d_i,\T\ov d_i\in\cm_\ell$.
But then the triangle $d_i\la\wh d_i\la\T\ov d_i$,
of the same $3\times3$ diagram, gives
that $\wh d_i,\T\wh d_i$ both
lie in $\cm_\ell*\cm_\ell\subset\cm_\ell$.
Lemma~\ref{L20.8}
allows us to conclude that, for
any pair of integers $m<m'$, both $>N$,
the triangle $c_m\la c_{m'}\la \wh d_{m,m'}$
is such that  $\wh d_{m,m'},\T\wh d_{m,m'}$
both lie in $\cm_\ell$.
This permits us to
conclude that
\be
\setcounter{enumi}{\value{enumiv}}
\item
In the sequence of triangles of (ii), the sequence $c_*$ is Cauchy.
\setcounter{enumiv}{\value{enumi}}
\ee
Which completes the proof of (ii).

As for (iii): we proved it as a step in
the proof of (ii).
\eprf

\rmk{R28.95}
Lemma~\ref{L28.19} produces examples of strong
triangles, and Remark~\ref{R28.905}
tells us that all strong
triangles are pre-triangles. Another source
of pre-triangles 
is mapping cones: given a morphism of pre-triangles
\[\xymatrix@C+20pt{
A\ar[r]^-{f}\ar[d]^{u} & B\ar[d]^{v}\ar[r]^-{g} &
C\ar[r]^-{h}\ar[d]^{w} & \T A\ar[d]^{\T u} \\
A'\ar[r]^-{f'} & B'\ar[r]^-{g'}    &
C'\ar[r]^-{h'}  & \T A'
}\]
then the mapping cone is also a pre-triangle
\[\xymatrix@C+40pt{
A'\oplus B
\ar[r]^-{\left(\begin{array}{cc}f' & v\\
0 & -g
\end{array}\right)} &
B'\oplus C\ar[r]^-{\left(\begin{array}{cc}g' & w\\
0 & -h
\end{array}\right)}    &
C'\oplus\T A\ar[r]^-{\left(\begin{array}{cc}h' & \T u\\
0 & -\T f
\end{array}\right)}
& \T A'\oplus\T B\ .
}\]
\ermk

In view of Remark~\ref{R28.95}, our next project is to learn
how to construct morphisms of pre-triangles. For this the next
little lemma is helpful.

\lem{L28.102}
Suppose we are given:
\be
\item
A homological object $B\in\fc(\cs)$. We remind the reader:
the fact that $B$ is homological means that $\Hom_{\MMod\cs}^{}(-,B)$
restricts to a homological functor on $\cs\op\subset\big[\MMod\cs\big]\op$.
\item
An object $A\in\fl(\cs)$. Assume also that we are given a
Cauchy sequence $a_*$ in $\cs$ with colimit $A$.
\ee
Then there exists an integer $n>0$ so that, for any integer $i\geq n$,
any map $a_i\la B$ factors
uniquely as $a_i\la A\la B$. More precisely: if we choose an integer $j>0$
with $B\in\cm_j^\perp$, then just choose $n$ to be an integer such that,
for all $i\geq n$,
the triangles $a_i\la a_{i+1}\la d_i$ have $\Tm d_i,d_i\in\cm_j^{}$.
\elem

\prf
Apply the homological functor $\Hom(-,B)$ to the triangles
$\Tm d_i\la a_i\la a_{i+1}\la d_i$. The hypotheses guarantee
that $\Hom(a_{i+1},B)\la\Hom(a_i,B)$ is an isomorphism whenever
$i\geq n$, allowing
us to extend any map $a_i\la B$, uniquely, to a map from $a_*$ to $B$.
\eprf

\cor{C28.105}
Assume we are given the following:
\be
\item
A strong triangle
$A\stackrel{f}\la B\stackrel{g}\la C\stackrel{h}\la \T A$
in the category $\fl(\cs)$,
which is the colimit of the image under $Y$
of the Cauchy sequence of triangles
$a_*\stackrel{f_*}\la b_*\stackrel{g_*}\la c_*\stackrel{h_*}\la \T a_*$
in the category $\cs$. 
\item
Composable morphisms
$A'\stackrel{f'}\la B'\stackrel{g'}\la C'\stackrel{h'}\la \T A'$
in the category $\fc(\cs)$, where the objects $A'$, $B'$ and $C'$ are all
homological.
\ee
Then there exists an integer $n>0$ such that, for any integer $i\geq n$,
any commutative
diagram
\[\xymatrix@C+20pt{\ar@{}[d]|{(1)}&
Y(a_i)\ar[r]^-{Y(f_i)}\ar[d]^{u_i} & Y(b_i)\ar[d]^{v_i}\ar[r]^-{Y(g_i)} &
Y(c_i)\ar[r]^-{Y(h_i)}\ar[d]^{w_i} & Y(\T a_i)\ar[d]^{\T u_i} &\\
&A'\ar[r]^-{f'} & B'\ar[r]^-{g'}    &
C'\ar[r]^-{h'}  & \T A' &
}\]
factors uniquely through a map 
\[\xymatrix@C+20pt{
A\ar[r]^-{f}\ar[d]^{u} & B\ar[d]^{v}\ar[r]^-{g} &
C\ar[r]^-{h}\ar[d]^{w} & \T A\ar[d]^{\T u} \\
A'\ar[r]^-{f'} & B'\ar[r]^-{g'}    &
C'\ar[r]^-{h'}  & \T A'
}\]

More precisely: if we choose $j>0$ so that $A',B',C',\T A'$ all belong to
$Y(\cm_j)^\perp$, it suffices to choose the integer $n>0$ large enough
so that, for any $i\geq n$, in the triangles
\[
a_i\stackrel{\alpha_i}\la a_{i+1}\la \ov d_i,\qquad
b_i\stackrel{\beta_i}\la b_{i+1}\la  d_i,\qquad
c_i\stackrel{\gamma_i}\la c_{i+1}\la \wh d_i
\]
we have
$\Tm \ov d_i,\ov d_i,\Tm d_i ,d_i,\Tm\wh d_i,\wh d_i$ belonging to
$\cm_j$.
\ecor

\prf
Lemma~\ref{L28.102} says that, with our choice of integer $n$, if we are
given an integer $i\geq n$ and a map from any of $a_i$, $b_i$ or $c_i$ to
any of $A'$, $B'$, $C'$ or $\T A'$, then the map factors uniquely
through the respective $a_i\la A$, $b_i\la B$ or $c_i\la C$. Applying this to
the maps $u_i:a_i\la A'$, $v_i:b_i\la B'$ and $w_i:c_i\la C'$ we factor
them uniquely as
\[
a_i\stackrel{\wt\alpha_i}\la A\stackrel u\la A'\,,\qquad
b_i\stackrel{\wt\beta_i}\la B\stackrel v\la B'\,,\qquad
c_i\stackrel{\wt\gamma_i}\la C\stackrel u\la C'\,,
\]
producing a
diagram
\[\xymatrix@C+20pt{
\ar@{}[dd]|{(2)}&a_i\ar[r]^-{f_i}\ar[d]^{\wt\alpha_i} & b_i\ar[d]^{\wt\beta_i}\ar[r]^-{g_i} &
c_i\ar[r]^-{h_i}\ar[d]^{\wt\gamma_i} & \T a_i\ar[d]^{\T\wt\alpha_i} &\\
 &A\ar[r]^-{f}\ar[d]^{u}& B \ar[r]^-{g}
\ar[d]^{v} &
C\ar[r]^-{h}\ar[d]^{w} & \T A\ar[d]^{\T u}&\\
&A'\ar[r]^-{f'} & B'\ar[r]^-{g'}   \ar[r] &
C'\ar[r]^-{h'}  & \T A'&
}\]
By construction we know that, when we delete the middle row of
(2), we are left
with diagram (1) in the statement of the Corollary---this
diagram commutes by hypothesis.  Deleting the bottom row of (2) leaves
us a commutative diagram, given by
the map from the triangle $a_i\la b_i\la c_i\la\T a_i$ to
the colimit of $a_*\la b_*\la c_*\la\T a_*$.
Thus the composites in each of the squares
\[\xymatrix@C+20pt{
 A\ar[r]^-{f}\ar[d]^{u}& B
\ar[d]^{v} &
B \ar[r]^-{g}
\ar[d]^{v} &
C\ar[d]^{w} & 
C\ar[r]^-{h}\ar[d]^{w} & \T A\ar[d]^{\T u}
\\
A'\ar[r]^-{f'} & B'&
B'\ar[r]^-{g'}   \ar[r] &
C' &
C'\ar[r]^-{h'}  & \T A'
}\]
give a pair of maps rendering equal the composites
\[\xymatrix@C+3pt{
a_i\ar[r]^-{\wt\alpha_i} & A\ar@<0.5ex>[r] \ar@<-0.5ex>[r] & B' & 
b_i\ar[r]^-{\wt\beta_i} & B\ar@<0.5ex>[r] \ar@<-0.5ex>[r] & C' & 
c_i\ar[r]^-{\wt\gamma_i} & C\ar@<0.5ex>[r] \ar@<-0.5ex>[r] & \T A'\ .  
}\]
The uniqueness assertion of Lemma~\ref{L28.102} gives that the three
squares must commute.
\eprf

\dfn{D28.110}
In the category $\fs(\cs)=\fl(\cs)\cap\fc(\cs)$, we declare a
sequence $A\stackrel{f}\la B\stackrel{g}\la C\stackrel{h}\la \T A$
to be a \emph{distinguished triangle} if
it is a strong triangle in $\fl(\cs)$
all of whose objects lie in $\fs(\cs)$.
\edfn

\rmk{R200976}
If $A\stackrel{f}\la B\stackrel{g}\la C\stackrel{h}\la \T A$
is a strong triangle in $\fl(\cs)$, in which $A$ and $B$
both belong to $\fs(\cs)$, then
$A\stackrel{f}\la B\stackrel{g}\la C\stackrel{h}\la \T A$
is a distinguished triangle in $\fs(\cs)$.

After all, Remark~\ref{R28.905} tells us that any
strong triangle is a pre-triangle, and
Remark~\ref{R28.109} tells us that if $A$ and
$B$ both lie in $\fc(\cs)$ then so does $C$.
Hence $C\in\fl(\cs)\cap\fc(\cs)=\fs(\cs)$, and
the strong triangle
$A\stackrel{f}\la B\stackrel{g}\la C\stackrel{h}\la \T A$
must be distinguished in $\fs(\cs)$.
\ermk

\thm{T28.128}
With the distinguished triangles as in Definition~\ref{D28.110},
the category $\fs(\cs)$ is triangulated. 
\ethm

\prf
We need to show that the axioms of triangulated categories are
satisfied. We begin with the obvious: for any object $A\in\fs(\cs)$ the
sequence $A\stackrel\id\la A\la 0\la\T A$ is clearly a triangle, just choose a
Cauchy sequence $a_*$ with colimit $A$ and consider the Cauchy sequence
of triangles $a_*\stackrel\id\la a_*\la0\la\T a_*$. Any isomorph of a triangle
is a triangle by definition.

Suppose we are given a morphism
$A\la B$ in $\fs(\cs)$. Lemma~\ref{L28.19} permits us to extend this
to a strong triangle $A\la B\la C\la \T A$ in the category $\fl(\cs)$, and 
Remark~\ref{R200976} shows that
this is a distinguished triangle in
$\fs(\cs)$.

This completes the proof of [TR1]. The axiom [TR2] is obvious,
the rotations of triangles in $\fs(\cs)$ are triangles.

It remains to prove [TR3] and [TR4'] as in
\cite[Definitions~1.1.2 and 1.3.13]{Neeman99}: we need to show that, given
a commutative diagram in the category $\fs(\cs)$ where
the rows are triangles
\[\xymatrix@C+20pt{
A\ar[r]^-{f}\ar[d]^{u} & B\ar[d]^{v}\ar[r]^-{g} &
C\ar[r]^-{h} & \T A \\
A'\ar[r]^-{f'} & B'\ar[r]^-{g'}    &
C'\ar[r]^-{h'}  & \T A'
}\]
we may complete it to a morphism of triangles, and even do so
in such a way that the mapping cone is a triangle.

OK: because the diagram lies in $\fs(\cs)\subset\fc(\cs)$ we may
choose an integer $j>0$ so that the objects
$A,B,C,\T A,\T B, A',B',C',\T A',\T B'$ all lie in $\cm_j^\perp$.
Next: because the rows are triangles we may choose Cauchy
sequences of triangles
$a_*\stackrel{f_*}\la b_*\stackrel{g_*}\la c_*\stackrel{h_*}\la \T a_*$
and
$a'_*\stackrel{f'_*}\la b'_*\stackrel{g'_*}\la c'_*\stackrel{h'_*}\la \T a'_*$
whose colimits are, respectively,
$A\stackrel{f}\la B\stackrel{g}\la C\stackrel{h}\la \T A$
and
$A'\stackrel{f'}\la B'\stackrel{g'}\la C'\stackrel{h'}\la \T A'$.
And since all these sequences are Cauchy we can, by ingoring
the early parts of the Cauchy sequences, assume that, with $x_*$ standing
for any of $a_*$, $b_*$, $c_*$, $a'_*$, $b'_*$ or $c'_*$, in any
triangle $x_m\la x_{m'}\la d$ the objects $\Tm d$, $d$ and $\T d$ lie
in $\cm_j$.

We are given the commutative diagram
\[\xymatrix@C+20pt{
a_1^{}\ar[r]^-{f_1}\ar[d]^{\wt\alpha_1^{}} & b_1^{}\ar[d]^{\wt\beta_1^{}}\\
A\ar[r]^-{f}\ar[d]^{u} & B\ar[d]^{v}\\
A'\ar[r]^-{f'} & B'
}\]
and, because $A'=\colim\,a'_*$ and $B'=\colim\, b'_*$, we may factor
the composite through some commutative diagram
\[\xymatrix@C+20pt{
a_1^{}\ar[r]^-{f_1}\ar[d]^{u_\ell^{}} & b_1^{}\ar[d]^{v_\ell^{}}\\
a'_\ell\ar[r]^-{f'_\ell}\ar[d]^{\wt\alpha'_\ell} & b'_\ell\ar[d]^{\wt\beta'_\ell}\\
A'\ar[r]^-{f'} & B'
}\]
In the triangulated category $\cs$ we may complete the
commutative diagram whose rows are triangles
\[\xymatrix@C+20pt{
a_1^{}\ar[r]^-{f_1}\ar[d]^{u_\ell^{}} & b_1^{}\ar[d]^{v_\ell^{}}\ar[r]^-{g_1^{}} &
c_1^{}\ar[r]^-{h_1} & \T a_1^{} \\
a'_\ell\ar[r]^-{f'_\ell} & b'_\ell\ar[r]^-{g'_\ell}    &
c'_\ell\ar[r]^-{h'_\ell}  & \T a'_\ell
}\]
to a morphism of triangles
\[\xymatrix@C+20pt{\ar@{}[d]|{(1)} &
a_1^{}\ar[r]^-{f_1}\ar[d]^{u_\ell^{}} & b_1^{}\ar[d]^{v_\ell^{}}\ar[r]^-{g_1^{}} &
c_1^{}\ar[r]^-{h_1}\ar[d]^{w_\ell^{}} & \T a_1^{}\ar[d]^{\T u_\ell^{}} &\\
& a'_\ell\ar[r]^-{f'_\ell} & b'_\ell\ar[r]^-{g'_\ell}    &
c'_\ell\ar[r]^-{h'_\ell}  & \T a'_\ell &
}\]
and even do so in such a way that the mapping cone is
a triangle. Now consider the composite
\[\xymatrix@C+20pt{
a_1^{}\ar[r]^-{f_1}\ar[d]^{u_\ell^{}} & b_1^{}\ar[d]^{v_\ell^{}}\ar[r]^-{g_1^{}} &
c_1^{}\ar[r]^-{h_1}\ar[d]^{w_\ell^{}} & \T a_1^{}\ar[d]^{\T u_\ell^{}} \\
a'_\ell\ar[r]^-{f'_\ell}\ar[d]^{\wt\alpha'_\ell}  &
b'_\ell\ar[r]^-{g'_\ell}\ar[d]^{\wt\beta'_\ell}     &
c'_\ell\ar[r]^-{h'_\ell}\ar[d]^{\wt\gamma'_\ell}
& \T a'_\ell\ar[d]^{\T\wt\alpha'_\ell} \\
A'\ar[r]^-{f'} & B'\ar[r]^-{g'}    &
C'\ar[r]^-{h'}  & \T A'
}\]
Corollary~\ref{C28.105} applies. Actually, we use the ``more precisely''
refinement, with $n=1$. The Corollary allows us to factor
the composite, uniquely, through
\[\xymatrix@C+20pt{\ar@{}[d]|{(2)} &
A\ar[r]^-{f}\ar[d]^{u} & B\ar[d]^{v}\ar[r]^-{g} &
C\ar[r]^-{h}\ar[d]^{w} & \T A\ar[d]^{\T u}& \\
&A'\ar[r]^-{f'} & B'\ar[r]^-{g'}    &
C'\ar[r]^-{h'}  & \T A'&
}\]
This already establishes [TR3], but we assert further that the mapping
cone is a triangle in $\fs(\cs)$.

To simplify the notation let us write the mapping cone of (2) as
$X\stackrel{\wt f}\la Y\stackrel{\wt g}\la Z\stackrel{\wt h}\la \T X$.
By Remark~\ref{R28.905} we know this to be a pre-triangle.
Note that, by our construction, the objects in
this pre-triangle, that is $X=A'\oplus B$,
$Y=B'\oplus C$, $Z=C'\oplus \T A$ and $\T X=\T A'\oplus \T B$,
all lie $\cm_j^\perp$.
We can furthermore express $X=A'\oplus B$ as $X\cong\colim\,(a'_*\oplus b_*)$
and $Y=B'\oplus C$ as $Y\cong\colim\,(b'_*\oplus c_*)$, that is we have
explicit Cauchy sequences converging to $X$ and $Y$. If we chop off the
sequences $a'_*$ and $b'_*$, deleting all the terms $a'_i$ and $b'_i$ with
$i<\ell$, we can even express $X=\colim\, x_*$ and $Y=\colim\,y'_*$
so that $x_1^{}=a'_\ell\oplus b_1^{}$ and $y'_1=b'_\ell\oplus c_1^{}$.
And the sequences $x_*$, $y'_*$ are such that, in the triangles
$x_m\la x_{m'}\la \ov d_{m,m'}$ and $y'_m\la y'_{m'}\la d_{m,m'}$, we have
\[
\ov d_{m,m'},\T\ov d_{m,m'},\Tm d_{m,m'}, d_{m,m'}
\qquad\text{ all belonging
  to }\qquad
\cm_j\ .
\]
And finally the morphism
from the mapping cone of (1) to the mapping cone of (2) rewrites as
\[\xymatrix@C+20pt{\ar@{}[d]|{(3)} &
x_1^{}\ar[r]^-{\wt f_1}\ar[d]^{\wt\alpha'_1} & y'_1\ar[d]^{\wt\beta'_1}\ar[r]^-{\wt g_1^{}} &
z_1^{}\ar[r]^-{\wt h_1}\ar[d]^{\wt\gamma'_1} & \T x_1^{}\ar[d]^{\T \wt\alpha'_1}& \\
& X\ar[r]^-{\wt f} & Y\ar[r]^-{\wt g}    &
Z\ar[r]^-{\wt h}  & \T X &
}\]
Where the top row is a distinguished
triangle in $\cs$, while the bottom row is a 
pre-triangle in $\fs(\cs)$.

Now we apply Lemma~\ref{L28.19} to the morphism $\wt f:X\la Y$ in $\fl(\cs)$.
Actually: we apply parts~(ii) and (iii) of Lemma~\ref{L28.19}.
We are given in $\cs$ a triangle
$x_1^{}\stackrel{\wt f_1}\la y'_1\stackrel{\wt g_1^{}}\la
z_1^{}\stackrel{\wt h_1}\la\T x_1^{}$, as well as a commutative square
\[\xymatrix@C+20pt{
x_1^{}\ar[r]^-{\wt f_1}\ar[d]^{\wt\alpha'_1} & y'_1\ar[d]^{\wt\beta'_1}\\
X\ar[r]^-{\wt f} & Y
}\]
and Cauchy sequences $x_*$, $y'_*$ with $X=\colim\,x_*$ and $Y=\colim\,y'_*$.
We may construct in $\cs$ a Cauchy sequence of triangles
$x_*\stackrel{\wt f_*}\la y_*\stackrel{\wt g_*}\la
z_*\stackrel{\wt h_*}\la\T x_*$,
extending the given triangle
$x_1^{}\stackrel{\wt f_1}\la y'_1\stackrel{\wt g_1^{}}\la
z_1^{}\stackrel{\wt h_1}\la\T x_1^{}$,
with $x_*$ the given sequence and $y_*$ a subsequence of $y'_*$
such that $y_1^{}=y'_1$,
and so that, in the Cauchy sequence $z_*$, we have
that the triangles $z_m\la z_{m'}\la\wh d_{m,m'}$ have
$\Tm\wh d_{m,m'},\wh d_{m,m'}$ both in $\cm_j$.
Let the colimit of
$x_*\stackrel{\wt f_*}\la y^{}_*\stackrel{\wt g_*}\la
z_*\stackrel{\wt h_*}\la\T x_*$
be written
$X\stackrel{\wt f}\la Y\stackrel{\wh g}\la\wh Z\stackrel{\wh h}\la \T X$.

Now we apply Corollary~\ref{C28.105} to the diagram in (3) and
the Cauchy sequence of triangles 
$x_*\stackrel{\wt f_*}\la y^{}_*\stackrel{\wt g_*}\la
z_*\stackrel{\wt h_*}\la\T x_*$. Actually: we apply the ``more precisely'' part
with $n=1$, to factor (3) uniquely through a morphism
\[\xymatrix@C+20pt{
X\ar[r]^-{\wt f}\ar@{=}[d] & Y\ar@{=}[d]\ar[r]^-{\wh g} &
\wh Z\ar[r]^-{\wh h}\ar[d]^{\ph} & \T X\ar@{=}[d]\ar[r]^{\T \wt f}
& \T Y\ar@{=}[d]\\
X\ar[r]^-{\wt f} & Y\ar[r]^-{\wt g}    &
Z\ar[r]^-{\wt h}  & \T X\ar[r]^{\T\wt f} &\T Y
}\]
The top row is a strong triangle in $\fl(\cs)$,
while the bottom row
is a pre-triangle in $\fc(\cs)$.
Thus both rows are exact in the abelian
category $\MMod\cs$. The 5-lemma
tells us that the map $\ph$ is an
isomorphism. Therefore $\wh Z\cong Z$ lies in $\fs(\cs)=\fl(\cs)\cap\fc(\cs)$,
the top  row is a triangle in $\fs(\cs)$, and
the bottom row, which is the mapping cone on the
morphism in (2), is isomorphic to the triangle in the top row.
\eprf

\rmk{R200978}
Let $\cs$ be an
essentially small triangulated category with a good metric.
The category $\MMod\cs$ is not usually triangulated, but it
contains the triangulated subcategories $\cs$ and $\fs(\cs)$.

More precisely, the triangulated category $\cs$ embeds
fully faithfully in $\MMod\cs$
via the Yoneda functor, and hence its essential image
$Y(\cs)\subset\MMod\cs$ has a triangulated structure. And the
triangulated structure on $\fs(\cs)$ comes from
Theorem~\ref{T28.128}. Of course, in the ambient
category $\MMod\cs$ we can form the intersection
\[
I(\cs)\eq Y(\cs)\cap\fs(\cs)\ ,
\]
and we assert that $I(\cs)$ is a triangulated subcategory,
both of $Y(\cs)$ and of $\fs(\cs)$.

Assume that $Y(a),Y(b)$ are objects of $I(\cs)$,
and $Y(f):Y(a)\la Y(b)$ is a morphism between them.
Of course, in the category $\cs$ we may complete to
a triangle $a\stackrel f\la b\stackrel g\la c\stackrel h\la \T a$, and then
form the Cauchy sequence of distinguished triangle
in $\cs$ where the connecting morphisms are identities
\[\xymatrix@C+20pt{
a\ar[r]^-{f}\ar[d]^{\id} & b\ar[d]^{\id}\ar[r]^-{g} &
c\ar[r]^-{h}\ar[d]^{\id} & \T a\ar[d]^{\id} &\\
a\ar[r]^-{f} & b\ar[r]^-{g}    &
c\ar[r]^-{h}  & \T a &
}\]
The colimit $Y(a)\la Y(b)\la Y(c)\la Y(\T a)$ is a strong
triangle in $\fl(\cs)$, where $Y(a)$ and $Y(b)$ are assumed
to lie in $\fs(\cs)\subset\fc(\cs)$. By
Remark~\ref{R200976}, it follows that
$Y(a)\la Y(b)\la Y(c)\la Y(\T a)$ is a distinguished
triangle in the category $\fs(\cs)$. In other
words: the completions of the morphism
$Y(f):Y(a)\la Y(b)$ in the category $I(\cs)$, to
distinguished triangles in the triangulated categories
$Y(\cs)$ and $\fs(\cs)$, coincide with each other.
\ermk

\section{In the presence of a good extension $\cs\la\ct$}
\label{S21}

In Sections~\ref{S20} and \ref{S28} we fixed a triangulated category
$\cs$ with a good metric, and out of it constructed and studied several
subcategories of $\MMod\cs$. But it turns out to be useful to embed
$\cs$ into other triangulated categories. Throughout, we
will assume given a fully faithful, triangulated functor
$F:\cs\la\ct$. Let us set up the conventions.

\ntn{N21.-100}
With $F:\cs\la\ct$ a fully faithful, triangulated functor, we let
$\cy:\ct\la\MMod\cs$ be the functor taking an object $A\in\ct$
to the functor $\cy(A)=\Hom\big(F(-),A\big)$. Clearly $\cy\circ F=Y$,
with $Y:\cs\la\MMod\cs$ as in the previous sections.
Because in this section we will be considering both the embedding
$Y:\cs\la\MMod\cs$ and the embedding $F:\cs\la\ct$, we will try to
be careful to not confuse $s\in\cs$ with its image under
either of
these embeddings.
\entn

We begin with the following. 

\obs{O21.-1}
Let $F:\cs\la\ct$ be as above. For every object $s\in\cs$,
Yoneda's lemma gives
\[
\Hom_{\MMod\cs}^{}\big(Y(s),\cy(t)\big)\eq\big[\cy(t)\big](s)
\eq\Hom_\ct^{}\big(F(s),t\big)\ ,
\]
which one can reformulate as saying that the natural map
\[\xymatrix{
\Hom_\ct^{}\big(F(s),t\big)\ar[r]^-\cy &
\Hom_{\MMod\cs}^{}\big(\cy F(s),\cy(t)\big)\ar@{=}[r] &
\Hom_{\MMod\cs}^{}\big(Y(s),\cy(t)\big)
}\]
is an isomorphism.
Specializing to  objects $s\in\cm_i\subset\cs$
gives
\[
\cy^{-1}\Big[Y(\cm_i)^\perp\Big]\eq\big[F(\cm_i)\big]^\perp\ .
\]
From Observation~\ref{O28.-1}, coupled with the fact
that $\cy^{-1}$ respects unions,
we obtain the formula
\[
\cy^{-1}\big(\fc(\cs)\big)\eq
\bigcup_{i\in\nn}^{}\,\,\big[F(\cm_i)\big]^\perp\ .
\]
From Observation~\ref{O28.-1} (i), (ii) and (iii) we
furthermore deduce
that $\cy^{-1}\big(\fc(\cs)\big)$ is a thick subcategory of $\ct$.
\eobs

\exm{E21.-3}
Let us specialize to the situation of Example~\ref{E20.3}(i):
the category $\ct$ has coproducts, there is 
a compact generator $H$ with $\Hom(H,\T^iH)=0$ for $i\gg0$,
and we are given a \tstr\ 
$\big(\ct^{\leq0},\ct^{\geq0}\big)$ in the preferred equivalence class.
As in Example~\ref{E20.3}(i) we set $\cs=\ct^c$ and
the good metric is given by $\cm_i=\ct^c\cap\ct^{\leq-i}$.

With $F:\ct^c\la\ct$ the natural embedding,
the inclusion $F(\cm_i)\subset\ct^{\leq-i}$ gives
$\ct^{\geq-i+1}=\big[\ct^{\leq-i}\big]^\perp\subset\big[F(\cm_i)\big]^\perp$. We
want to prove an inclusion in the other direction.
Note that, because $H$ is compact and
the \tstr\ is in the preferred equivalence class, there
is an integer $n>0$ with $\T^nH\in\ct^{\leq0}$, and hence
$\genul H{}{-n-i}\subset\ct^c\cap\ct^{\leq-i}=\cm_i$.
Therefore
$\big[F(\cm_i)\big]^\perp\subset\big[\genul H{}{-n-i}\big]^\perp=\big[\ogenul H{}{-n-i}\big]^\perp$.
The definition of the \tstr\
$\big(\ct_H^{\leq0},\ct_H^{\geq0}\big)$ generated by $H$ is
that $\ct_H^{\leq0}=\ogenul H{}0$, and as the \tstr\
$\big(\ct^{\leq0},\ct^{\geq0}\big)$  
is in the preferred equivalence class it is equivalent
to $\big(\ct_H^{\leq0},\ct_H^{\geq0}\big)$, and there is an
integer $n'>0$ with
$\ogenul H{}{-n-i}=\ct_H^{\leq-n-i}\supset\ct^{\leq-n-n'-i}$. Taking
perpendiculars, and combining with
the earlier inclusions, we deduce
\[
\ct^{\geq-i+1}\subset \big[F(\cm_i)\big]^\perp\subset\big[\ogenul H{}{-n-i}\big]^\perp\subset\big[\ct^{\leq-n-n'-i}\big]^\perp=\ct^{\geq-n-n'-i+1}\ .
\]
Taking the union over $i\in\nn$ we discover
$\cup_{i\in\nn}^{}\big[F(\cm_i)\big]^\perp=\ct^+$.
Appealing to
the formula for $\cy^{-1}\big(\fc(\cs)\big)$ given in
Observation~\ref{O21.-1},
this combines to
\[
\cy^{-1}\big(\fc(\cs)\big)\eq\ct^+\ .
\]
\eexm

\exm{E21.-5}
Still with $\ct$ being a triangulated category with coproducts and with
$\cs=\ct^c$, we can consider a different good metric. More explicitly, with
$H$ still a compact generator for $\ct$ we can let
$\cm_i=\big[\genu H{}{-i,i}\big]^\perp$.

First we note that this is---up to equivalence---the good metric studied in
Krause~\cite{Krause18}. Krause formulates his theory
terms of Cauchy sequences and not metrics,
hence let us sketch the translation.
If a sequence $E_*$
is Cauchy with respect to the  good metric above then in the triangle
$E_n\la E_{n+1}\la D_n$  we have $\Tm D_n,D_n\in\cm_i$ for $n\gg0$.
This implies that, for all
objects $G\in{\genu H{}{-i,i}}$, the functor $\Hom(G,-)$ takes the maps
$E_n\la E_{n+1}$ to isomorphisms whenever $n\gg0$.
Since every $G\in\ct^c=\gen H{}=\cup_{i=1}^\infty\genu H{}{-i,i}$
belongs to some $\genu H{}{-i,i}$, it follows that a Cauchy sequence
with respect to the good metric is Cauchy in Krause's sense.

Conversely: if the sequence is Cauchy in Krause's sense then, with
$\ov G=\oplus_{\ell=-i-1}^i\T^\ell H$ we have that $\Hom(\ov G,-)$ takes
 $E_n\la E_{n+1}$ to isomorphisms for all $n\gg0$. Put
 $G=\oplus_{\ell=-i}^i\T^\ell H$; since both $\Tm G$ and $G$ are
 direct summands of $\ov G$, we have that both $\Hom(G,-)$ and
$\Hom(\Tm G,-)$ take $E_n\la E_{n+1}$
 to isomorphisms, or to rephrase this, $\Hom(G,-)$ takes both $E_n\la
 E_{n+1}$ and $\T E_n\la\T E_{n+1}$ to isomorphisms. We conclude
 that, in the triangle 
$E_n\la
 E_{n+1}\la D_n$, we have $D_n\in G^\perp$ for $n\gg0$. As
$G^\perp=\big[\oplus_{\ell=-i}^i\T^\ell H\big]^\perp=\big[\genu H{}{-i,i}\big]^\perp=\cm_i$
 the sequence is Cauchy in the good metric.

We can, of course, try to compute what the theory developed here
yields when applied to the good metric that underlies Krause's
construction. The only
case I have computed in detail is
when $\ct=\D(\Mod R)$ with $R$ a noetherian
ring. We can choose the compact generator to
be $R\in\D(\Mod R)$, and then the categories $\cm_i$ come down to
\[
\cm_i\eq\{X\in\D^b(\proj R)\mid H^\ell(X)=0\text{ whenever }-i\leq\ell\leq i\}\ .
\]
Because $\T^jR\in\cm_i$ whenever $|j|>i$, we have
that any object $X\in\big[F(\cm_i)\big]^\perp$ must have
$H^j(X)=0$ if $|j|>i$, in other words $X\in\D(\Mod R)^{\geq-i}\cap\D(\Mod R)^{\leq i}$.
But the
category $\cm_i$ also contains good approximations
for every object of the form $\T^{-i-1}M$, where
$M\in\mod R$. Precisely: choose a resolution for $\T^{-i-1}M$ by
finitely generated, projective $R$--modules,
that is a complex
\[\xymatrix{
\cdots P^{i-1}\ar[r] &P^{i}\ar[r] &
P^{i+1}\ar[r] & 0\ar[r] &\cdots
}\]
whose only cohomology is $M$ in degree $i+1$.
Then form the brutal truncation, deleting everything in degree $<-i-1$.
We obtain
an object in $P^*\in\D^b(\proj R)$  with only two nonvanishing
cohomology groups,
$H^{i+1}(P^*)=M$ and $H^{-i-1}(P^*)$. Hence
$P^*\in\cm_i$.  The triangle
$(P^*)^{\leq-i-1}\la P^*\la \T^{-i-1}M\la \T(P^*)^{\leq-i-1}$ tells us that,
for $X\in\big[F(\cm_i)\big]^\perp\subset\D(\Mod R)^{\geq-i}$ the map
$\Hom(\T^{-i-1}M,X)\la\Hom(P^*,X)=0$ must be an isomorphism.
But $X\in\big[F(\cm_i)\big]^\perp$ also belongs to $\D(\Mod R)^{\leq i}$,
and the vanishing of
$\Hom(\T^{-i-1}M,X)$
for every finite $R$--module $M$ guarantees 
that $X$ must be isomorphic to a bounded complex
of injective $R$--modules, vanishing outside degrees $-i\leq j\leq i$.

Now use  the formula of Observation~\ref{O21.-1} for
$\cy^{-1}\big(\fc(\cs)\big)$. It tells us that, for
Krause's good metric
on $\cs=\D^b(\proj R)$, the category $\cy^{-1}\big(\fc(\cs)\big)\subset\D(\Mod R)$
turns out to be the category of bounded complexes of injectives.
\eexm

So far we have
computed a couple of examples, we have worked out what the category
$\cy^{-1}\big(\fc(\cs)\big)$ comes down to in the special
cases of Example~\ref{E20.3}(i) and of Krause's good metric on $\D^b(\proj
R)$. 
We would like to also say something
about $\fl(\cs)$, and for this it helps to restrict the class of
embeddings $\cs\la\ct$ we consider. This leads us to
the following.

\dfn{D21.1}
Let $F:\cs\la\ct$ be a fully faithful triangulated functor
between triangulated categories.
We say that $F:\cs\la\ct$ is a \emph{good extension
with respect to the metric} if
\be
\item
In $\ct$ we can form the homotopy colimit of
$F(E_*)$, for any Cauchy sequence
$E_*$ in the category $\cs$. We will
spell this out more fully in
Explanation~\ref{E21.19178}.
\item
For any Cauchy sequence $E_*$ in
$\cs$, the natural map $\colim\, Y(E_i)\la\cy\big(\hoco F(E_i)\big)$
is an isomorphism.
\ee
\edfn

\xpl{E21.19178}
We should elaborate on the somewhat cryptic
hypothesis made in  
Definition~\ref{D21.1}(i).

Definition~\ref{D21.1}(i) states that,
in the category $\ct$, we can
form the homotopy colimit, of
the image under $F$ of any Cauchy sequence
$E_*$ in the category $\cs$.
This can be expanded into the 
following two assertions:
\be
\item
Given any Cauchy sequence $E_*$, in the category
$\cs$, the
countable coproduct $\coprod_{i=1}^\infty F(E_i)$
exists in $\ct$.
\item
The homotopy colimit $\hoco F(E_*)$ is defined,
up to non-canonical isomorphism, by the
triangle
\[\xymatrix{
\ds\coprod_{i=1}^\infty F(E_i)\ar[rr]^-{\id-\sh} & &
\ds\coprod_{i=1}^\infty F(E_i)\ar[r] &\hoco F(E_*)
\ar[r] &
\ds\coprod_{i=1}^\infty \T F(E_i)\ .
}\]
\ee
\expl

\dis{D21.19180}
In both Definition~\ref{D21.1} and
Explanation~\ref{E21.19178}
we were careful \emph{not to assume} that
the category $\ct$ has arbitrary small coproducts.

In the examples of this article this care is
unnecessary. The generality in which
Definition~\ref{D21.1} was
formulated will play a role only
in
subsequent work.

Since we have chosen not to assume that $\ct$
has coproducts, we need to tell the reader
how to modify the one standard construction
we will want to use. Recall: if $B\subset\ct$ is
any collection of objects, the full subcategory
$\loc B$ in normally defined to be the smallest
localizing subcategory containing $B$, meaning
the smallest triangulated subcategory of $\ct$,
containing $B$, and closed under coproducts.

Since we are not assuming that $\ct$ has all
small coproducts, the modified definition of $\loc B$
is that it is minimal subject
to the hypotheses that it contains $B$, is
a triangulated subcategory of $\ct$, is
closed in $\ct$ under those coproducts that
exist in $\ct$, and is also closed in $\ct$ under
direct summands

If all small coproducts exist in
  $\ct$, then
  closure under direct summands follows automatically
  from the other assumptions.
\edis

We will have occasion to use the following
little lemma more than once.

\lem{L21.19182}
Let $\ct$ be a triangulated category. For
any triangulated subcategory $\cs\subset\ct$,
if $X\in\cs^\perp$ then $X\in{\loc\cs}^\perp$.
\elem

\prf
As $\cs$ is a triangulated subcategory we know that
$\T^n\cs=\cs$ for all integers $n$.
But $X\in\cs^\perp=\big(\T^n\cs\big)^\perp$ implies
$\T^{-n}X\in\cs^\perp$, and hence $\cs^\perp$ contains
$X(-\infty,\infty)$, with the conventions of
\cite[Notation~1.11]{Neeman17}. Thus
$\cs\subset{^\perp\big[X(-\infty,\infty)\big]}$.
But as ${^\perp\big[X(-\infty,\infty)\big]}$ is
a triangulated subcategory of $\ct$, containing $\cs$,
closed under
those coproducts that exist in $\ct$, and
closed under
direct summands in $\ct$, we conclude that
$\loc\cs\subset{^\perp\big[X(-\infty,\infty)\big]}$.
\eprf

\exm{E21.3}
Let $\ct$ be a triangulated category, and let
$\ct^c$ be the full subcategory of compact
objects in $\ct$. By this we mean: an object $X\in\ct$
belongs to $\ct^c$ if and only if $\Hom_\ct^{}(X,-)$
commutes with those coproducts that exist in $\ct$.
If $\cs\subset\ct^c$ is any
triangulated subcategory, then
the embedding $F:\cs\la\ct$ is a good extension
for any good metric on $\cs$
such that, for Cauchy sequences $E_*$
in $\cs$,  the coproduct $\coprod_{i\in\nn}F(E_i)$
exists in $\ct$. For example: if $\ct$ contains
all countable coproducts, then
there is no restriction on the
metric.

The fact that this is a good extension
follows from \cite[Lemma~2.8]{Neeman96};
the proof shows that, for any sequence
$E_*$ in $\ct$ whose
homotopy colimit exists as
in Explanation~\ref{E21.19178}, the map
$\colim\,\cy(E_*)\la\cy(\hoco E_*)$ is
an isomorphism.
\eexm

\exm{E21.5}
Now let $\ct$ be a weakly approximable triangulated category,
and choose a \tstr\ 
$\big(\ct^{\leq0},\ct^{\geq0}\big)$ in the preferred equivalence class.
As in Example~\ref{E20.3}(ii) let $\cs=\big[\ct^b_c\big]\op$,
and the good metric is given by $\cm_i\op=\ct^b_c\cap\ct^{\leq-i}$.

I assert that the embedding $F:\big[\ct^b_c\big]\op\la\ct\op$ is
a good extension---the reader can find the proof in
\cite[Lemmas~~3.1]{Neeman24}.
\eexm

The definition of $\fl(\cs)$ was to take colimits in
the category $\MMod\cs$, but since good extensions
$F:\cs\la\ct$ are such that $\ct$ is assumed to contain
certain homotopy colimits, the following
makes sense.

\dfn{D21.7}
Let $F:\cs\la\ct$ be a good extension.
The full subcategory $\fl'(\cs)\subset\ct$ has for objects all
the isomorphs in $\ct$ of homotopy
colimits of the images under $F$ of Cauchy sequences in $\cs$.
\edfn

\obs{O21.9}
If $F:\cs\la\ct$ is  a good extension, then 
the functor $\cy:\ct\la\MMod\cs$ restricts, on objects, to an essential
surjection $\text{Ob}\big(\fl'(\cs)\big)\la\text{Ob}\big(\fl(\cs)\big)$.
After all: by Definition~\ref{D21.7} the functor $\cy$
takes homotopy colimits of Cauchy sequences to colimits.
\eobs

\exm{E21.10}
We return to Example~\ref{E20.3}(i), where
$\cs=\ct^c$ and the good metric is given by $\cm_i=\ct^c\cap\ct^{\leq-i}$.
We assert that
in this case the category $\fl'(\cs)$ turns out to be $\ct^-_c$.
The reader can find this 
in \cite[Lemma~8.5]{Neeman24}.
\eexm

\exm{E21.10.5}
With $\ct$ a triangulated category with coproducts and
a single compact generator $H$, we can let $\cs=\ct^c$ be as above,
but endow
$\cs$ with Krause's good metric---see Example~\ref{E21.-5}.
By Example~\ref{E21.3} the pair $\cs\subset\ct$ is a good
extension. We can form the category $\fl'(\cs)$, but I have only
computed it when $\ct=\D(\Mod R)$ for a noetherian ring $R$.

Every object $E\in\cs=\D^b(\proj R)$ has bounded cohomology, with
$H^j(E)$ a finite $R$--module for every $i$. In any Cauchy
sequence, with respect to Krause's good metric, the cohomology
eventually stabilizes. Therefore for any $X\in\cl'(\cs)$
and any $j\in\zz$ we have that $H^i(X)$ is a finite $R$--module.
In symbols: $\cl'(\cs)\subset\D_{\mod R}^{}(R)$,
the category of all complexes of $R$--modules with finite
cohomology modules.

I assert that this inclusion is an equality.
Suppose $X$ belongs to $\D_{\mod R}^{}(R)$, I want to exhibit $X$ as
the homotopy colimit of a Cauchy sequence. To this end pick an integer $i>0$
and 
consider the map $X^{\leq i}\la X$ from the truncation with
respect to the standard \tstr\ on $\D(\Mod R)$. The object $X^{\leq i}$ is
bounded above and has finite cohomology modules, hence admits a
resolution by finitely generated projectives---there is in
$\D(\Mod R)$ an isomorphism $P\la X^{\leq i}$, with $P\in\D^-(\proj R)$.
Now take the brutal truncation, killing everything
in degree $<-i-1$ to obtain a map $E_i\la X^{\leq i} \la X$ with
$E_i\in\D^b(\proj R)$. The
functor $H^j(-)$ takes this map to an isomorphism whenever $-i\leq j\leq i$,
and these assemble to a Cauchy sequence with homotopy colimit $X$.
\eexm

Now that we have seen a few examples, let us turn to general results that
hold under the assumptions introduced. 

\lem{L21.11}
Let $F:\cs\la\ct$ be a good extension, as in Definition~\ref{D21.1}.
Suppose $E$ belongs to the category $\fl'(\cs)$ of
Definition~\ref{D21.7}, and let $X\in\ct$ be arbitrary.
Then the natural map $\Hom(E,X)\la\Hom\big(\cy(E),\cy(X)\big)$
extends to a short exact sequence
\[\xymatrix{
0\ar[r] & K(E,X) \ar[r] &
\Hom(E,X)\ar[r] &\Hom\big(\cy(E),\cy(X)\big)\ar[r] & 0\ .
}\]
Moreover: if $E_*$ is a Cauchy sequence with $E\cong\hoco F(E_*)$,
then there is an isomorphism
$K(E,X)\cong\clim^1\,\Hom\big(\T F(E_i),X\big)$.
\elem

\prf
Choose a Cauchy sequence
$E_*$ in $\cs$ with $E=\hoco F(E_*)$, and consider
the commutative square
\[\xymatrix@C+40pt{
\Hom(E,X)\ar[r]\ar[d] & \clim\Hom(F(E_i),X)\ar[d]^{\wr}\\
\Hom\big(\cy(E),\cy(X)\big)\ar[r]^-\sim & \clim\Hom\big(Y(E_i),\cy(X)\big) \ .
}\]
The vertical map on the right is an isomorphism
by Observation~\ref{O21.-1}, and the bottom
horizontal map is an isomorphism because 
\begin{eqnarray*}
\Hom\big(\cy(E),\cy(X)\big) &=& \Hom\big[\cy\big(\hoco F(E_i)\big)\,,\,\cy(X)\big]\\
&=& \Hom\big[\colim\,Y(E_i)\,,\,\cy(X)\big]\\
&=& \clim\Hom\big[Y(E_i)\,,\,\cy(X)\big]
\end{eqnarray*}
where the second isomorphism is because we're dealing with
a good extension. We are therefore reduced to showing
that the top horizontal map
is a surjection, and computing its kernel.

For this recall the definition of homotopy colimits: the homotopy
colimit $E=\hoco F(E_*)$ sits in a triangle
\[\xymatrix{
\ds\coprod_{i>0}F(E_i) \ar[rr]^-{1-\text{shift}} &&
\ds\coprod_{i>0}F(E_i) \ar[r] &
E \ar[r] & 
\ds\coprod_{i>0}\T F(E_i) 
}\]
and, applying the functor $\Hom(-,X)$, we obtain a
short exact sequence
\[\xymatrix{
0\ar[r] & \clim^1\,\Hom\big(\T F(E_i),X\big) \ar[r] &
\Hom(E,X)\ar[r] &\clim\,\Hom\big(F(E_i),X\big)\ar[r] & 0
}\]
This completes the proof of the Lemma.
\eprf

\rmk{R21.900}
The next few results will specialize Lemma~\ref{L21.11}
to the case where $X$ also belongs to $\fl'(\cs)$. In other
words we will study the essential surjection
$\cy:\fl'(\cs)\la\fl(\cs)$
of Observation~\ref{O21.9}.
\ermk

\cor{C21.902}
The essential surjection $\cy:\fl'(\cs)\la\fl(\cs)$
of Observation~\ref{O21.9}
is full.
\ecor

\prf
Fullness says that the map
$\Hom_{\fl'(\cs)}^{}(E,E')\la\Hom_{\fl(\cs)}^{}\big(\cy(E),\cy(E')\big)$
is surjective, and this follows immediately by
specializing Lemma~\ref{L21.11} to the case where $E'$
also belongs to $\fl'(\cs)$.
\eprf

\dfn{D21.904}
A morphism $f:A\la B$ in the category $\fl'(\cs)$ is
called \emph{$\cs$--phantom} if $\cy(\ph)=0$.
\edfn

\lem{L21.906}
The collection of $\cs$--phantom maps in $\fl'(\cs)$ forms
a square-zero two-sided ideal.
\elem

\nin
We remind the reader: being a two-sided ideal
means that the $\cs$--phantom maps are closed under
addition, as well as composition, on the left
or on the right, by arbitrary
morphisms in $\fl'(\cs)$. And being square-zero
means that the composite of any two $\cs$--phantom
maps vanishes.

\prf
The fact that the $\cs$--phantom maps
form a two-sided ideal is immediate from their being
the kernel of the
functor $\cy:\fl'(\cs)\la\fl(\cs)$.

Now for the ideal being square-zero.
Let $f:A\la B$ be an $\cs$--phantom map.
Because $A$ belongs to $\fl'(\cs)$, there exists in
$\ct$ a triangle forming the top row of the
diagram below
\[\xymatrix@C+10pt{
\ds\coprod_{i=1}^\infty F(a_i)\ar[rr]^-{\id-\sh} &&
\ds\coprod_{i=1}^\infty F(a_i)\ar[r]\ar[dr]_-0 &
A\ar[r]\ar[d]^f&
\ds\coprod_{i=1}^\infty \T F(a_i)\\
& & &B &
}\]
The map $f:A\la B$ is such that $\cy(f)=0$, meaning
that $\Hom\big(F(\cs),-\big)$ annihilates $f$. Therefore the slanted
composite in the diagram above vanishes, and the map $f$
must factor as
\[\xymatrix@C+30pt{
A\ar[r]&
\ds\coprod_{i=1}^\infty \T F(a_i)\ar[r] &
B
}\]
If $g:B\la C$ is another phantom map, then the composite
$A\stackrel f\la B\stackrel g\la C$ factors through
the vanishing composite
\[\xymatrix@C+30pt{
\ds\coprod_{i=1}^\infty \T F(a_i)\ar[r] &
B \ar[r]^-g & C
}\]
\eprf

\cor{C21.908}
The functor $\cy:\fl'(\cs)\la\fl(\cs)$ reflects isomorphisms.
And since it is both essentially surjective and full, it follows
that on objects it gives a bijection between the isomorphism
classes of objects in $\fl'(\cs)$ and the
isomorphism classes of objects in $\fl(\cs)$.
\ecor

\prf
Suppose $f:A\la B$ is a morphism in $\fl'(\cs)$, such that $\cy(f)$ is an
isomorphism. Because $\cy$ is full, we can choose a morphism
$g:B\la A$ such that $\cy(g)=\cy(f)^{-1}$. Therefore $gf:A\la A$ and $fg:B\la B$ are  morphisms
with $\cy(gf)=\id$ and with $\cy(fg)=\id$.

But then $h=\id-gf$ and $k=\id-fg$ are morphisms $h:A\la A$
and $k:B\la B$, with both $\cy(h)=0$ and $\cy(k)=0$. That is
both $h$ and $k$ are $\cs$--phantom maps.
By Lemma~\ref{L21.906} we have that
$h^2=0=k^2$. This makes $gf=\id-h$ and $fg=\id-k$
both isomorphisms, with inverses $(\id+h)$ and $(\id+k)$
respectively. 
Therefore $f:A\la B$ has both a right and a left inverse,
and must be an isomorphism in $\fl'(\cs)$.

As for the other statements in the Corollary: if $A$ and
$B$ are objects in $\fl'(\cs)$, with $\cy(A)$ and $\cy(B)$
isomorphic, then we can choose an isomorphism
$\ph:\cy(A)\la\cy(B)$. As the functor $\cy$ is full, we can
choose a morphism $f:A\la B$ with $\cy(f)=\ph$ an
isomorphism. But then $f$ is an isomorphism, proving the
assertion about the bijection of isomorphism classes.
\eprf

For now this is all we want to say about the relationship
between the categories $\fl'(\cs)$ and $\fl(\cs)$. We move
on to the following.

\lem{L21.13}
If $E$ is an object of $\fl'(\cs)$ and $X\in\ct$ is an object
with $\cy(X)\in\cc(\cs)$,
then
the natural map $\Hom(E,X)\la\Hom\big(\cy(E),\cy(X)\big)$ is an
isomorphism.
\elem

\prf
Choose a Cauchy sequence $E_*$ in the category
$\cs$ with $E=\hoco\, F(E_*)$.
Since $X\in\ct$ is such that
$\cy(X)\in\fc(\cs)$, Observation~\ref{O21.-1}
allows us to choose
an integer $n$ with
$\Hom_\ct^{}\big(F(\cm_n),X\big)=0$.
Because the sequence $E_*$ is Cauchy there exists an integer $M>0$
so that, for all $M\leq m<m'$, the triangle
$D_{m,m'}\la\T E_m\la \T E_{m'}\la \T D_{m,m'}$ in the category $\cs$ has
$D_{m,m'},\T D_{m,m'}\in\cm_n$. Applying the exact functor $\Hom\big(F(-),X\big)$
to this triangle we have that the map
$\Hom\big(F(\T E_{m'}),X\big)\la\Hom\big(F(\T E_{m}),X\big)$ is an isomorphism
whenever $M\leq m<m'$, and hence $\clim^1\,\Hom\big(F(\T E_i),X\big)=0$.
The current Lemma now follows Lemma~\ref{L21.11}.
\eprf

\cor{C21.15}
The restriction of $\cy:\ct\la\MMod\cs$, to the subcategory
$\fl'(\cs)\cap\cy^{-1}\big(\fc(\cs)\big)\subset\ct$,
induces an equivalence with the category $\fs(\cs)$
of Definition~\ref{D20.11}(iii).
\ecor

\prf
In Observation~\ref{O21.9} we noted that the functor $\cy$
yields an essential surjection
$\text{Ob}\big(\fl'(\cs)\big)\la\text{Ob}\big(\fl(\cs)\big)$,
and restricting to the inverse image
of $\fs(\cs)=\fl(\cs)\cap\fc(\cs)$
will yield an essential surjection
$\text{Ob}\big[\fl'(\cs)\cap\cy^{-1}\big(\fc(\cs)\big)\big]\la\text{Ob}\big(\fs(\cs)\big)$. So on objects the functor is essentially surjective.

On the other hand Lemma~\ref{L21.13} tells us that, for $E,X$ in
$\fl'(\cs)\cap\cy^{-1}\big(\fc(\cs)\big)$, the map
$\Hom(E,X)\la\Hom\big(\cy(E),\cy(X)\big)$ is an
isomorphism. Thus $\cy$ is fully faithful on the subcategory.
\eprf

And now we come to the point.

\thm{T20.17}
Let $\cs$ be a triangulated category with a good metric $\{\cm_i\}$, and
let $F:\cs\la\ct$ be a good extension.
Then the category $\fl'(\cs)\cap\cy^{-1}\big(\fc(\cs)\big)$ is a
triangulated subcategory of $\ct$, and the natural map
\[\xymatrix@C+30pt{
\cy\,\,:\,\,\fl'(\cs)\cap\cy^{-1}\big(\fc(\cs)\big)\ar[r] &\fs(\cs)
}\]
is a triangulated equivalence, where the category on the left has
the triangulated structure in inherits from being a triangulated
subcategory of $\ct$, and the category $\fs(\cs)$ has the
triangulated structure of Definition~\ref{D28.110}.
\ethm

\prf
The fact that $\cy$ is an equivalence of categories was proved
in Corollary~\ref{C21.15}, we only have to worry about the
triangulated structure. 
Let $A$, $B$ be objects in $\fl'(\cs)\cap\cy^{-1}\big(\fc(\cs)\big)$,
and suppose 
$A\stackrel{f}\la B\stackrel{g}\la C\stackrel{h}\la \T A$
is a triangle in $\ct$. We need to prove
that $C$ belongs to $\fl'(\cs)\cap\cy^{-1}\big(\fc(\cs)\big)$,
and that
$\cy(A)\stackrel{\cy(f)}\la \cy(B)\stackrel{\cy(g)}\la \cy(C)
\stackrel{\cy(h)}\la \cy(\T A)$
is a triangle in $\fs(\cs)$.

Remark~\ref{R28.109} guarantees that, in the pre-triangle
$\cy(A)\stackrel{\cy(f)}\la \cy(B)\stackrel{\cy(g)}\la \cy(C)
\stackrel{\cy(h)}\la \cy(\T A)$, the object
$\cy(C)$ must belong to $\fc(\cs)$. As the objects $\cy(A)$,
$\cy(B)$, $\cy(C)$ and $\cy(\T A)$ all belong
to $\fc(\cs)$, we choose and fix an integer $n$ so that all
four objects belong to $\big[Y(\cm_n)\big]^\perp$. Because
the objects $\cy(A)$ and $\cy(B)$ belong to
$\fl(\cs)$ we may choose Cauchy sequences converging to
them: that is we pick a Cauchy sequence $a_*$ with
$\cy(A)=\colim\,Y(a_*)$ and a Cauchy sequence $b'_*$ with
$\cy(B)=\colim\,Y(b'_*)$. Passing to subsequences if necessary,
we make sure that for all integers $0<m<m'$, in the
triangles $a_m\la a_{m'}\la \ov d_{m,m'}$ and
$b_m\la b_{m'}\la d_{m,m'}$ we have that
$\ov d_{m,m'},\T\ov d_{m,m'},\Tm d_{m,m'},d_{m,m'}$ all
belong to $\cm_n$.

In the category $\fl(\cs)$ we have the
map $\cy(f):\cy(A)\la\cy(B)$, to which we may apply
Lemma~\ref{L28.19}: we may choose a subsequence
$b_*$ of $b'_*$, a sequence of triangles
$a_*\stackrel{f_*}\la b_*\stackrel{g_*}\la c_*\stackrel{h_*}\la \T a_*$
in $\cs$, do it in such a way that in the
triangles $c_m\la c_{m'}\la \wh d_{m,m'}$
we have $\Tm\wh d_{m,m'},\wh d_{m,m'}\in \cm_n$,
and ensure that the colimit of 
$Y(a_*)\stackrel{Y(f_*)}\la Y(b_*)\stackrel{Y(g_*)}\la Y(c_*)\stackrel{Y(h_*)}\la Y(\T a_*)$
 is a strong triangle 
$\cy(A)\stackrel{\cy(f)}\la \cy(B)\stackrel{\wt g}\la \wt C
\stackrel{\wt h}\la \cy(\T A)$.
 
In particular the construction gives us a commutative square
in $\MMod\cs$
 \[\xymatrix@C+30pt{
Y(a_1) \ar[r]^-{Y(f_1)}\ar[d] & Y(b_1)\ar[d] \\
\cy(A)\ar[r]^-{\cy(f)} & \cy(B)
}\]
which must be the image under $\cy$ of a commutative square in $\ct$
\[\xymatrix@C+30pt{
F(a_1^{}) \ar[r]^-{F(f_1)}\ar[d] & F(b_1)\ar[d] \\
A\ar[r]^-{f} & B
}\]
This last square may be extended to a morphism
of triangles in $\ct$
\[\xymatrix@C+30pt{
F(a_1^{}) \ar[r]^-{F(f_1)}\ar[d] & F(b_1)\ar[d]\ar[r]^-{F(g_1)}\ar[d] & F(c^{}_1)\ar[d]\ar[r]^-{F(h_1)}\ar[d] & F(\T a^{}_1)\ar[d] \\
A\ar[r]^-{f} & B \ar[r]^-{g} & C\ar[r]^-{h} & \T A
}\]
and applying the functor $\cy$ we deduce a
commutative diagram in $\MMod\cs$
\[\xymatrix@C+30pt{
Y(a_1^{}) \ar[r]^-{Y(f_1)}\ar[d] & Y(b_1)\ar[d]\ar[r]^-{Y(g_1)} & Y(c^{}_1)\ar[d]\ar[r]^-{Y(h_1)} & Y(\T a^{}_1)\ar[d] \\
\cy(A)\ar[r]^-{\cy(f)} & \cy(B) \ar[r]^-{\cy(g)} & \cy(C)\ar[r]^-{\cy(h)} & \cy(\T A)
}\]
And now we apply Corollary~\ref{C28.105}, or rather
the ``more precisely'' version with $n=1$,
to factor this map, uniquely,
through
\[\xymatrix@C+30pt{
\cy(A) \ar[r]^-{\cy(f)}\ar@{=}[d] & \cy(B)\ar@{=}[d]\ar[r]^-{\wt g} & \wt C\ar[d]^\ph \ar[r]^-{\wt h}\ar[d] & \cy(\T A)\ar@{=}[d] \ar[r]^-{\cy(\T f)} & \cy(\T B)\ar@{=}[d]\\
\cy(A)\ar[r]^-{\cy(f)} & \cy(B) \ar[r]^-{\cy(g)} & \cy(C)\ar[r]^-{\cy(h)} & \cy(\T A)\ar[r]^-{\cy(\T f)} & \cy(\T B)
}\]
The 5-lemma forces $\ph$ to be an isomorphism. The top row is
a triangle in $\fs(\cs)$ by construction, and the isomorphism
tells us that so is the bottom row.

It remains to prove that $\fl'(\cs)\cap\cy^{-1}\big(\fc(\cs)\big)$ is
a triangulated subcategory of $\ct$,
concretely we still need to check that $C$ belongs to
$\fl'(\cs)\cap\cy^{-1}\big(\fc(\cs)\big)$. What we know so far is that
$\cy(C)\in\fc(\cs)$, or equivalently that 
$C\in
\cy^{-1}\big(\fc(\cs)\big)$; 
it remains to prove that
$C\in\cl'(\cs)$. Because $F:\cs\la\ct$
is a good extension we have isomorphisms
$\wt C=\colim\,Y(c_*)\cong\cy(\hoco c_*)$.
This makes $\ph$ a morphism
$\ph:\cy(\hoco c_*)\la\cy(C)$ with $\hoco c_*\in\fl'(\cs)$ and
$\cy(C)\in\fc(\cs)$.
Lemma~\ref{L21.13} allows us to lift the isomorphism $\ph$ to a (unique)
morphism $\rho:\hoco c_*\la C$ in the category $\ct$. And since
$\cy(\rho)=\ph$ is an isomorphism, it follows that in the triangle
$D\la \hoco c_*\la C$ we have $\cy(D)=0$, or
to rephrase we have $D\in\cs^\perp$.
By Lemma~\ref{L21.19182} it follows that
$D\in{\loc\cs}^\perp$, with
$\loc\cs$ as in Discussion~\ref{D21.19180}.

On the other hand we have that $A$, $B$ and
$\hoco c_*$ belong to $\cl'(\cs)\subset\loc\cs$,
the triangle $A\la B\la C\la\T A$ informs us that $C\in\loc\cs$,
and from the triangle $D\la \hoco c_*\la C$ we learn that $D\in\loc\cs$.
The identity map $\id:D\la D$ is a morphism from $D\in{\loc\cs}$
to $D\in{\loc\cs}^\perp$ and must vanish.
Hence $D=0$, the map $\rho:\hoco c_*\la C$
is an isomorphism, and this exhibits $C$ as isomorphic to
$\hoco c_*\in\fl'(\cs)$.
\eprf

\rmk{R29876}
As in Remark~\ref{R200978}, we can wonder about the
category $I(\cs)=Y(\cs)\cap\fs(\cs)$.

Suppose we are given a good extension $F:\cs\la\ct$.
In  Notation~\ref{N21.-100} we noted that $\cy\circ F=Y$,
and in Corollary~\ref{C21.15} we proved that
$\cy$ gives an equivalence of
$\cl'(\cs)\cap\cy^{-1}\big(\fs(\cs)\big)$ with $\fs(\cs)$.
This means that, up to an equivalence induced
by the functor $\cy$, the category $I(\cs)$ is given
by the formula
\[
\begin{array}{ccccc}
I(\cs)&\eq& Y(\cs)\cap\fs(\cs)&\quad\cong\quad&
F(\cs)\cap\cl'(\cs)\cap\cy^{-1}\big(\fc(\cs)\big)\\
& & &\eq&
F(\cs)\cap\cy^{-1}\big(\fc(\cs)\big)
\end{array}
\]
where the last equality is because $F(\cs)\subset\cl'(\cs)$.
In Remark~\ref{R200978} we saw that $I(\cs)$ is a
triangulated subcategory of both $\cs$ and
$\fs(\cs)$.
By Theorem~\ref{T20.17},
the triangulated
structure on $I(\cs)\subset\fs(\cs)$
therefore agrees with the triangulated
structure on
\[
F(\cs)\cap\cy^{-1}\big(\fc(\cs)\big)\sub
\fl'(\cs)\cap\cy^{-1}\big(\fc(\cs)\big)\sub\ct\ ,
\]
allowing us to compute.
\ermk

\section{The example of $\ct^c$ determining $\ct^b_c$}
\label{S22}

It's time to see what the generality of Sections~\ref{S20}, \ref{S28}
and \ref{S21} reduces to
in some concrete examples. But first a few global conventions for
this section.

\ntn{N22.1}
Throughout this section $\ct$ will
be a triangulated category with coproducts,
we will assume there is a compact generator $H\in\ct$ with
$\Hom(H,\T^iH)=0$ for $i\gg0$, and we will suppose given a \tstr\
$\big(\ct^{\leq0},\ct^{\geq0}\big)$ in the preferred equivalence class.
\entn

\exm{E22.3}
Let the conventions be as in Notation~\ref{N22.1}.
In Example~\ref{E20.3}(i) we studied the following: we put
$\cs=\ct^c$, and let the good metric be given by $\{\cm_i=\ct^c\cap\ct^{\leq-i}\}$.
In Example~\ref{E21.3} we learned that the embedding $F:\ct^c\la\ct$
is a good extension, Example~\ref{E21.-3} teaches us that
the category  $\cy^{-1}\big(\fc(\cs)\big)$ turns out to be
$\ct^{+}$, while in Example~\ref{E21.10} we saw that $\fl'(\cs)=\ct^-_c$.
This means that 
\[
\fl'(\cs)\cap\cy^{-1}\big(\fc(\cs)\big)\eq
\ct^-_c\cap\ct^+\eq\ct^b_c\ .
\]
Theorem~\ref{T20.17} tells us that the functor $\cy$ induces a
a triangulated equivalence of
$\ct^b_c=\fl'(\cs)\cap\cy^{-1}\big(\fc(\cs)\big)$
with the category $\fs(\cs)$. This computes for us what $\fs(\cs)$ turns
out to be, as a triangulated category, as long as the good metric is
as in  Example~\ref{E20.3}(i).
\eexm

From Example~\ref{E22.3}
we learn that, in the generality given in Notation~\ref{N22.1},
the triangulated category $\ct^b_c$ is fully determined by the
category $\ct^c$ \emph{together with its good metric.} The way we defined
the good metric was to use the embedding into $\ct$; our definition
was $\cm_i=\ct^c\cap\ct^{\leq-i}$. While it's true that, up to
equivalence of good metrics, the \tstr\ doesn't matter much---equivalent
{\it t}--structures induce equivalent good metrics, and the
{\it t}--structures in
the preferred equivalence class are all equivalent---the preferred
equivalence class of {\it t}--structures is defined on $\ct$, not
on $\ct^c$.

Hence the reader might wonder if there is some way to define the good metric
on $\ct^c$ without referring to the embedding into $\ct$. 
We begin with the following.

\rmd{R22.99}
A \emph{classical generator} of a
triangulated category $\cs$ is
an object $G\in\cs$ with $\cs=\gen G{}$. From
\cite[Lemma~0.13(ii)]{Neeman24} it follows that, given two classical
generators $G$ and $H$, there exists an integer $A$ with
$H\in\genu G{}{-A,A}$ and $G\in\genu H{}{-A,A}$.
\ermd

\lem{L22.101}
Suppose $\ct$ is a compactly generated triangulated category.
Then any classical generator
for $\ct^c$ is a compact generator for $\ct$.
\elem

\prf
Let $G\in\ct^c$ be a classical generator. Then $\ct^c=\gen G{}$, and
$\ct=\loc{\ct^c}=\ogen G$.
\eprf

\dfn{D22.103}
Let $\cs$ be a triangulated category, and assume $G\in\cs$ is a
classical generator. We define two good metrics $\{\cl_i\}$, $\{\cn_i\}$
by the formulas
\be
\item
$\cl_i=\genul G{}{-i}$.
\item
$\cn_i=\big[\genuf G{}{-i}\big]^\perp$.
\ee
\edfn

\rmk{R22.105}
The Cauchy sequences with respect to the good metric $\{\cn_i\}$ manifestly
agree with those of Remark~\ref{R0.5}.

In this generality all that's clear is that, up to equivalence, the
good metrics $\{\cl_i\}$ and $\{\cn_i\}$
don't depend on the choice of classical generator. This
follows immediately from Reminder~\ref{R22.99}, more specifically
from the fact that, given two classical generators $G$ and $H$,
we can find an integer $A>0$ with
$H\in\genu G{}{-A,A}$ and $G\in\genu H{}{-A,A}$.
\ermk

\rmk{R22.105.5}
Now return to the situation of 
Notation~\ref{N22.1}. Put $\cs=\ct^c$ and consider the good metrics
$\{\cl_i\}$, $\{\cn_i\}$ of Definition~\ref{D22.103}, as well as
the good metric $\{\cm_i\}$ of Example~\ref{E20.3}(i).
The classical generator $G\in\cs=\ct^c$
 of Definition~\ref{D22.103} is a compact
object in $\ct$, and the \tstr\
$\big(\ct^{\leq0},\ct^{\geq0}\big)$
of Notation~\ref{N22.1} 
belongs to the preferred equivalence class. By
\cite[Observation~0.16 and Lemma~2.8]{Neeman24} there
exists an integer $A$ with $\T^A G\in\ct^{\leq0}$ and with 
$\Hom\big(\T^{-A}G,\ct^{\leq0}\big)=0$. We deduce the
following.
\be
\item
$\qquad\cl_{i+A}\eq\genul G{}{-i-A}\sub\ct^c\cap\ct^{-i}\eq \cm_i$
\item
 $\qquad \cm_{i+A}\eq\ct^c\cap\ct^{-i-A}\sub\ct^c\cap\big[\genuf G{}{-i}\big]^\perp 
  \eq\cn_i$
\ee
By Lemma~\ref{L22.101} the classical generator $G\in\ct^c$ is a compact
generator of $\ct$. 
But then the fact that the \tstr\ $\big(\ct^{\leq0},\ct^{\geq0}\big)$
is in the preferred equivalence class says that {\it t}--structures
$\big(\ct^{\leq0},\ct^{\geq0}\big)$ and 
$\big(\ct_G^{\leq0},\ct_G^{\geq0}\big)$
must be equivalent, which gives an integer $B>0$
for which the inclusion
$\ct^{\leq-i-B}\subset\ct_G^{\leq-i}=\ogenul G{}{-i}$ holds. Now
 \cite[Proposition~1.9]{Neeman17} gives the second last equality in
\[
\cm_{i+B}\eq\ct^c\cap\ct^{\leq-i-B}\sub\ct^c\cap\ogenul G{}{-i}\eq
\genul G{}{-i}\eq\cl_{i}\ .
\]  
That is, without assuming any approximability we have
$\{\cl_i\}\cong\{\cm_i\}$ and $\{\cm_i\}\preceq\{\cn_i\}$.
\ermk

If we assume weak approximability we can do better.

\pro{P22.109}
Assume $\ct$ is a weakly approximable triangulated category, and let
$\cs=\ct^c$. Then the good metric $\{\cn_i\}$
on $\cs$ given in
Definition~\ref{D22.103}
is equivalent to
the good metric $\{\cm_i\}$ of Example~\ref{E20.3}(i).
\epro

\prf
In Remark~\ref{R22.105.5} we noted the inequalities
$\{\cm_i\}\preceq\{\cn_i\}$; what needs proof is
the reverse inequality.
By Lemma~\ref{L22.101} the classical generator $G\in\ct^c$ is a compact
generator of $\ct$.
 By
\cite[Lemma~3.9(iv)]{Burke-Neeman-Pauwels18} there is an integer $B>0$ with
$\big[\genuf G{}{-i}\big]^\perp\subset\ct^{\leq-i+B}$, hence
\[
\cn_{i+B}\eq\ct^c\cap\big[\genuf G{}{-i-B}\big]^\perp\sub\ct^c\cap\ct^{\leq-i}\eq
\cm_i\ .
\]
\eprf

\rmk{R22.856}
Let us remain with the conventions of this section, that is $\cs=\ct^c$ where
$\ct$ satisfies the assumptions of Notation~\ref{N22.1}: it's natural to ask
what happens with Krause's good metric. As I have
said before, I have only worked out completely what
happens in the case $\ct=\D(\Mod R)$ with $R$ a noetherian ring.

In Example~\ref{E21.-5} we saw that, with respect to the good extension
$F:\D^b(\proj R)\la\D(\Mod R)$, the subcategory $\cy^{-1}\big(\fc(\cs)\big)$
turns out to be $\D^b(\Inj R)$, the category of bounded complexes
of injective modules. In Example~\ref{E21.10.5} we learned
that the category $\fl'(\cs)$ is $\D_{\mod R}^{}(R)$,
the full subcategory of $\D(\Mod R)$ of complexes whose cohomology
modules are finite. This makes
the category $\fl'(\cs)\cap\cy^{-1}\big(\fc(\cs)\big)$ equal to
$\D^b_{\mod R}(\Inj R)$, meaning the objects
are the bounded complexes of injectives
whose cohomology modules are finite.

Assume $R$ is commutative and
there is a dualizing complex $C\in\D^b(\mod R)$, meaning
a complex with a bounded injective resolution, and such that
$\Hom(-,C)$ induces an equivalence $\D^b(\mod R)\op\la\D^b(\mod R)$.
Then $\Hom(-,C)$ takes an object of $\fl'(\cs)\cap\cy^{-1}\big(\fc(\cs)\big)$
to a complex in $\D^b(\mod R)$ with a bounded flat resolution, hence
a bounded projective resolution---in other words
to an object of $\D^b(\proj R)$. Thus the
category $\fs(\cs)\cong\fl'(\cs)\cap\cy^{-1}\big(\fc(\cs)\big)$
is equivalent to $\D^b(\proj R)\op\cong\D^b(\proj R)$.

That is: if $R$ is a commutative noetherian ring with
a dualizing complex and if $\cs=\D^b(\proj R)$,
and if the \emph{good metric we use is Krause's metric,}
then $\fs(\cs)=\cs$. Whereas with the good metric
we have been using in the current paper, we have
\[
\fs(\cs)\eq\fs\big(\D^b(\proj R)\big)\eq\D^b(\mod R)\ .
\]
This illustrates the philosophical point we made
in Remark~\ref{R0.212121}: the choice of metric
really does matter. And for an elegant, enhancement-free
passage from $\D^b(\proj R)$ to $\D^b(\mod R)$,
the metric of this article is much preferable
to Krause's metric.
\ermk

\rmk{R309183}
In Remark~\ref{R200978} we introduced the category
$I(\cs)=Y(\cs)\cap\fs(\cs)$ and showed that it
is a triangulated subcategory, both of $Y(\cs)\cong\cs$
and of $\fs(\cs)$. And in Remark~\ref{R29876}
we learned that, given a good extension $F:\cs\la\ct$,
there is a simple formula for $I(\cs)$ as
a triangulated subcategory of $\ct$.

Let us use the formula where $\cs=\D(\proj R)$
with the metric of this article---that is either
of the equivalent metrics $\{\cl_i\}$ and $\{\cn_i\}$
of Definition~\ref{D22.103}---and the good extension
is the natural functor $F:\D^b(\proj R)\la\D(\Mod R)$.
Then the
formula yields
\[
\begin{array}{ccccl}
I(\cs)&\eq& \cs\cap\ct^+&\eq&\D^b(\proj R)\cap\D^+(\Mod R)\\
 & & & \eq & \D^b(\proj R)\ .
\end{array}
\]
Now assume $R$ is a commutative, noetherian
ring with a dualizing complex, we give the
category $\D^b(\proj R)$ Krause's metric,
and the good extension is
still $F:\D^b(\proj R)\la\D(\Mod R)$. Then the formula
of Remark~\ref{R29876} yields
\[
I\big(\D^b(\proj R)\big)\eq 
\D^b(\proj R)\cap\D^b(\Inj R)\ .
\]
Unless $R$ is a Gorenstein ring, the object $R$ belongs
to $\D^b(\proj R)$ but does not have a bounded injective
resolution, that is
\[
R\quad\in\quad\D^b(\proj R)-I\big(\D^b(\proj R)\big)\ .
\]
That is the inclusion
\[
I\big(\D^b(\proj R)\big)\quad\subsetneq\quad\D^b(\proj R) 
\]
is proper.

Now note: by Remark~\ref{R22.856} we know that,
with $\cs=\D^b(\proj R)$ being given
Krause's metric, 
there is an equivalence of categories
\[
\fs\big(\D^b(\proj R)\big)\quad\cong\quad\D^b(\proj R)\ .
\]
However: both are subcategories of $\MMod\cs$, and
the embeddings must be different---after all the
intersection $I(\cs)$ is a proper subcategory of $\cs$.
\ermk

\section{Coherent approximable categories, and passing from $\ct^b_c$
back to $\ct^c$}
\label{S29}

The reader may have noticed that in the treatment so far we have been strangely
reticent about Example~\ref{E20.3}(ii). Recall: under the
hypotheses of Notation~\ref{N22.1} we can look at the category
$\cs=\big[\ct^b_c\big]\op$, endow it with the good metric
$\cm_i\op=\ct^b_c\cap\ct^{\leq-i}$, and consider the embedding
$F:\cs=\big[\ct^b_c\big]\op\la\ct\op$. In Example~\ref{E21.5} we
mentioned that, as long
as $\ct$ is weakly approximable,
this is a good extension with respect to the good metric. But since
then there has been silence---no mention of the example. The reason
is that without further hypotheses there isn't much to say, what we can
prove is that $\fl'(\cs)$ is contained in $\big[\ct^-_c\big]\op$; this
follows from \cite[Lemma~3.3]{Neeman24}.

The reason we can't say much more
without hypotheses is simple: without some restriction, I see
no reason for the category $\ct^b_c$ to be nonzero---we have
proved that $\fl'(\cs)\subset\ct^-_c$, but to expect
an inclusion in the other direction there better be some
nonzero objects in $\ct^b_c$. And by way of
cautionary example: if
$R$ is a DGA, and $H^i(R)=0$ for $i>0$, then
\cite[Remark~4.3]{Neeman24} teaches us
that the category
$\ct=\ch^0\big(\Mod R\big)$ is an approximable
triangulated category---in this example we still have a reasonable
understanding of the preferred equivalence
class of \tstr{s}. But then
Bondarko and
Vostokov~\cite[Corollary~4.3]{Bondarko-Vostokov20}
improved on this and showed that,
\emph{assuming only that $H^i(R)=0$ for $i>1$,}
the category $\ct=\ch^0\big(\Mod R\big)$
is weakly approximable.
And in this situation the preferred equivalence
class of \tstr{s} is mysterious, we have no
idea how to compute $\ct^b_c$, and in particular
don't know if it is nonzero.

\dfn{D29.1}
Suppose $\ct$ is a triangulated category with coproducts, and assume
it has a compact generator $H$ with $\Hom(H,\T^iH)=0$ for $i\gg0$.
We declare $\ct$ to be \emph{coherent} if there exists an
integer $N>0$, and a \tstr\
$\big(\ct^{\leq0},\ct^{\geq0}\big)$ in the preferred equivalence class,
such that: for every object $X\in\ct^-_c$ there exists a triangle
$A\la X\la B$ with $A\in\ct^-_c\cap\ct^{\leq0}$ and $B\in\ct^-_c\cap\ct^{\geq-N}=\ct^b_c\cap\ct^{\geq-N}$.
\edfn

\rmk{R29.3}
The definition is clearly robust. If one \tstr\
has an integer $N$ as above, then so
does any equivalent \tstr. The integer will of course depend on the
\tstr.
\ermk

\exm{E29.7}
Let $\ct$ be a weakly approximable triangulated
category, 
and assume there exists, in
the preferred equivalence class of
\tstr{s} on $\ct$, a \tstr\ which restricts
to a \tstr\ on $\ct^-_c$. That is: we have
in the preferred equivalence class
a \tstr\ $(\ct^{\leq0},\ct^{\geq0})$ such that,
if $F$ is an object of $\ct^-_c$ and
$F^{\leq0}\la F\la F^{\geq1}$ is the triangle
that comes from truncating with respect
to the \tstr\
$(\ct^{\leq0},\ct^{\geq0})$, then
$F^{\leq0}$ and $F^{\geq1}$ both belong to
$\ct^-_c$.

Then of course the category $\ct$ is coherent.
For the \tstr\ $(\ct^{\leq0},\ct^{\geq0})$ of
the paragraph above, which is assumed to belong
to the preferred equivalence class, we
can choose the integer $N>0$
of Definition~\ref{D29.1} to be $N=1$.
\eexm

\exm{E29.5}
Let $X$ be a quasicompact, quasiseparated,
coherent scheme\footnote{We remind the reader:
a scheme $X$ is \emph{coherent} if
the structure sheaf $\co_X^{}$ is a
coherent sheaf of rings. In particular: all
noetherian schemes are coherent.}.
Let $Z\subset X$ be a Zariski-closed subset
with quasicompact complement.
From \cite[Theorem~3.2(iv)]{Neeman22A}
we learn that the category $\ct=\Dqcs Z(X)$ is
weakly approximable, while
\cite[Theorem~3.2(iii)]{Neeman22A}
informs us that the standard \tstr\
is in the preferred equivalence class.
And \cite[Theorem~3.3]{Neeman22A}
amounts to a computation of $\ct^-_c$.
Of course \cite[Theorem~3.3]{Neeman22A}
is stated in gorgeous generality, there
is no assumption that the scheme $X$ is
coherent. In the special case,
where $X$ is assumed coherent, the
computation simplifies to
$\ct^-_c=\dmcohs Z(X)$. That is: $\ct^-_c$
is the derived category whose objects $F$
are complexes of sheaves of $\co_X^{}$ modules,
where the restriction of $F$ to the open set $X-Z$ is
acyclic, and whose cohomology sheaves $\ch^i(F)$
are coherent,
and vanish when $i\gg0$.

Since the standard \tstr\ on $\ct$ is in the preferred
equivalence class, and obviously restricts to
a \tstr\ on $\ct^-_c$, we are in the situation of
Example~\ref{E29.7}. The weakly
approximable triangulated category $\Dqcs Z(X)$ is
coherent.
\eexm

\exm{E29098793}
Let $R$ be a left-coherent
ring---possibly noncommutative.
The fact that $\ct=\D(\Mod R)$ is approximable and
the computation that $\ct^-_c=\D^-(\mod R)$ may be
found in \cite[Example~4.1]{Neeman24},
and it is also observed that the standard \tstr\
is in the preferred equivalence class.
Thus Example~\ref{E29.7} also applies here to tell
us that $\D(\Mod R)$ is a coherent approximable category.
\eexm

\rmk{E29098793.567}
There exist rings $R$, which are not left-coherent,
but whose derived category $\D(\Mod R)$ is coherent.
The examples I know come from Rickard
(unpublished): it is possible
to construct a pair of derived-equivalent rings $R$
and $S$, with $S$ left-coherent but $R$ not. The
derived category
$\D(\Mod R)\cong\D(\Mod S)$ is a coherent triangulated
category because of the left-coherence of $S$,
see Example~\ref{E29098793}.
\ermk

\exm{E29098794}
An analysis similar to that
of Examples~\ref{E29.7},
\ref{E29.5} and
\ref{E29098793}, beginning with
\cite[Example~4.2]{Neeman24}, shows that
the homotopy category
of spectra $\ct$ is coherent and approximable.
\eexm

\lem{L29.50509}
Let $\ct$ be a coherent, weakly approximable triangulated category,
and choose a \tstr\  $\big(\ct^{\leq0},\ct^{\geq0}\big)$
in the preferred equivalence class.
For any object $X\in\ct^-_c\cap\ct^{\leq0}$ there exists a Cauchy
sequence $B_*$ 
in $\ct^b_c\cap\ct^{\leq0}$ with $X=\holim B_*$.
\elem

\prf
Let $N$ be the integer whose existence is guaranteed by
Definition~\ref{D29.1}, and suppose we are given
an object $X\in\ct^-_c\cap\ct^{\leq0}$. Definition~\ref{D29.1}
permits us to produce, for every
integer $i>0$, a triangle $A_i\la X\la B_i$ with
$A_i\in\ct^-_c\cap\ct^{\leq-i(N+1)}$ and $B_i\in\ct^-_c\cap\ct^{\geq-i(N+1)-N}$.
Therefore the composite $A_{i+1}\la X\la B_i$ is a morphism
from $A_{i+1}\in\ct^{\leq-i(N+1)-N-1}$ to $B_i\in\ct^{\geq-i(N+1)-N}$
and must vanish. We can factor $X\la B_i$ through $X\la B_{i+1}\la B_i$,
creating an inverse system in $\ct^b_c$. Because
$A_i\in\ct^{\leq-i(N+1)}$ the functor $(-)^{\geq-i(N+1)+1}$
takes the maps $X\la B_i$ to an isomorphism, showing that
$B_i\in\ct^b_c\cap\ct^{\leq0}$ and
that the sequence $B_*$ is Cauchy. And
\cite[Proposition~3.2]{Neeman24} shows that $X\la\holim B_*$
is an isomorphism.
\eprf

\pro{P29.9}
Let $\ct$ be a coherent, weakly approximable triangulated category,
and choose a \tstr\  $\big(\ct^{\leq0},\ct^{\geq0}\big)$
in the preferred equivalence class.
With $\cs=\big[\ct^b_c\big]\op$, with the good metric
$\cm_i\op=\ct^b_c\cap\ct^{\leq-i}$, and with the good extension
$F:\cs=\big[\ct^b_c\big]\op\la\ct\op$, we
have the following.
\be
\item
$\fl'(\cs)=\big[\ct^-_c\big]\op$.
\item
$\fl'(\cs)\cap\cy^{-1}\big(\fc(\cs)\big)=\big[\ct^c\big]\op$
\ee
\epro

\prf
The inclusion $\fl'(\cs)\subset\big[\ct^-_c\big]\op$ is contained in
\cite[Lemma~3.3]{Neeman24},
and the inclusion $\big[\ct^-_c\big]\op \subset \fl'(\cs)$ follows
from
Lemma~\ref{L29.50509}. This proves (i).

Now for (ii). Suppose $X\in\ct^-_c$ belongs to $\Big[\big[F(\cm_i)\big]^\perp\Big]\op$,
that is $\Hom(X,B)=0$ for all $B\in\cm_i\op=\ct^b_c\cap\ct^{\leq-i}$.
Choose any object $Y\in\ct^-_c\cap\ct^{\leq-i-1}$. By Lemma~\ref{L29.50509}
we can express $Y$ as
$Y=\holim B_n$, with $B_n\in\ct^b_c\cap\ct^{\leq-i-1}$. 
Hence $Y$ sits in a triangle 
\[\xymatrix@C+40pt{
\ds\prod_{n=1}^\infty \Tm B_n\ar[r] & Y\ar[r] &\ds\prod_{n=1}^\infty  B_n
}\]
Since $\Tm B_n,B_n$ belong to $\ct^b_c\cap\ct^{\leq-i}=\cm_i\op$ for every $n$,
the functor $\Hom(X,-)$ annihilates both outside terms in
the triangle above. Therefore $\Hom(X,Y)=0$.
We learn that, if $X\in\ct^-_c$  belongs to $\Big[\big[F(\cm_i)\big]^\perp\Big]\op$,
then $\Hom(X,Y)=0$ for all $Y\in\ct^-_c\cap\ct^{\leq-i-1}$. But
the fact that $X\in\ct^-_c$ means that there must
exist a triangle $E\la X\la Y$ with $E\in\ct^c$ and
$Y\in\ct^-_c\cap\ct^{\leq-i-1}$. The vanishing of the map $X\la Y$ forces
$X$ to be a direct summand of $E\in\ct^c$. This proves
the inclusion
$\fl'(\cs)\cap\cy^{-1}\big(\fc(\cs)\big)\subset\big[\ct^c\big]\op$.

Now for the reverse inclusion. We have that
$\big[\ct^c\big]\op\subset\big[\ct^-_c\big]\op=\fl'(\cs)$,
hence what needs proof is that $\big[\ct^c\big]\op$ belongs
to $\Big[\big[F(\cm_i)\big]^\perp\Big]\op$ for some $i\in\nn$.
But this is easy, any object $X\in\ct^c$ satisfies
$\Hom\big(X,\ct^{\leq-i}\big)=0$ for some $i>0$,
and $F(\cm_i)\subset\ct^{\leq-i}$.
\eprf

\section{The category $\ct^b_c$ determines $\ct^c$}
\label{S1}

In Proposition~\ref{P29.9} we saw that, if $\ct$ is coherent
and weakly approximable, then
the category $\ct^b_c$ \emph{together with its good metric} determines
$\ct^c$. We want to find a recipe to produce the good metric
just from the triangulated category $\ct^b_c$.

\ntn{N1.-10}
In this section it is quite easy to become confused by perpendiculars.
The convention will be: if $P\subset\ct$ is any subset, then
$P^\perp$ means the perpendicular of $P$ in $\ct$.
If it so happens that $P\subset\ct^b_c\subset\ct$, our notation
for the perpendicular of $P$ in $\ct^b_c$ will be $\ct^b_c\cap P^\perp$.
\entn

\lem{L1.1}
Let $\ct$ be a coherent triangulated category,
and choose a \tstr\  $\big(\ct^{\leq0},\ct^{\geq0}\big)$ in the
preferred equivalence class.
For any object $H\in\ct^-_c$ and any integer $m\geq0$ we can find
an object $G\in\ct^b_c$ with
$H^\perp\supset \ct^{\geq-m}\cap G^\perp$.
\elem

\prf
The coherent hypothesis permits us to construct a triangle
$A\la H\la G$, with $A\in\ct^{\leq-m-1}$ and
$G\in\ct^-_c\cap\ct^{\geq-m-1-N}\subset\ct^b_c$.
The result now follows because
$H^\perp\supset A^\perp\cap G^\perp\supset\ct^{\geq-m}\cap G^\perp$.
\eprf

\lem{L1.99}
Let $\ct$ be a coherent triangulated category.
In the partial order on the subcategories of $\ct^b_c$, introduced in
Definition~\ref{D0.13}, the subcategories
$\ct^b_c\cap\big[\genul G{}{0}\big]^\perp$ with $G\in\ct^b_c$
have (up to equivalence) a
unique minimal member $\cq(\ct^b_c)$.
The subcategory $\cq(\ct^b_c)$ is equivalent, in the partial order,
to $\ct^b_c\cap\ct^{\geq0}$, where $\big(\ct^{\leq0},\ct^{\geq0}\big)$
is a \tstr\ in the preferred equivalence class.
\elem

\prf
To produce a $G$ which gives (up to equivalence)
the unique minimal
$\ct^b_c\cap\big[\genul G{}{0}\big]^\perp$,
we apply Lemma~\ref{L1.1} to a compact generator $H\in\ct$, to the
\tstr\ $\big(\ct_H^{\leq0},\ct_H^{\geq0}\big)$, and to the integer $m=0$.
We learn that we may find an object
$G\in\ct^b_c$ with $H^\perp\supset \ct_H^{\geq0}\cap G^\perp$.
Suspending $i$ times, with $i\geq 0$,
gives
\[
(\T^iH)^\perp\quad\supset\quad \ct_H^{\geq-i}\cap(\T^iG)^\perp
\quad\supset\quad \ct_H^{\geq0}\cap(\T^iG)^\perp
\]
Intersecting over $i\geq0$ we obtain
\[
\bigcap_{i=0}^\infty (\T^iH)^\perp \quad\supset\quad
\ct_H^{\geq0}\cap\left[\bigcap_{i=0}^\infty \big[(\T^iG)^\perp\right]
\]
which rewrites as
$\ct_H^{\geq1}\supset \ct_H^{\geq0}\cap\big[\genul G{}{0}\big]^\perp$.
Suspending $i+1$, times, with $i\geq0$, gives
$\ct_H^{\geq-i}\supset\ct_H^{\geq-i-1}\cap\big[\genul G{}{-i-1}\big]^\perp$,
and induction allows us to prove, for any $i\geq0$,
$\ct_H^{\geq1}\supset\ct_H^{\geq-i}\cap\big[\genul G{}{0}\big]^\perp$. After all
it's true for $i=0$, and the inductive step follows from
\begin{eqnarray*}
\ct_H^{\geq1} &\supset &\ct_H^{\geq-i}\cap\big[\genul G{}{0}\big]^\perp\\
&\supset &\ct_H^{\geq-i-1}\cap\big[\genul G{}{-i-1}\big]^\perp\cap \big[\genul G{}{0}\big]^\perp
\\
&=&\ct_H^{\geq-i-1}\cap\big[\genul G{}{0}\big]^\perp
\end{eqnarray*}
Now taking the union over $i>0$ gives
$\ct_H^{\geq1}\supset\ct^+\cap\big[\genul G{}{0}\big]^\perp$.
Intersecting with $\ct^b_c$ we have
$\ct^b_c\cap\ct_H^{\geq1}\supset\ct^b_c\cap\big[\genul G{}{0}\big]^\perp$.
This gives us that, in the partial order of
Definition~\ref{D0.13}, we have
$\ct^b_c\cap\big[\genul G{}{0}\big]^\perp\preceq
\ct^b_c\cap\ct_H^{\geq0}$.

On the other hand for any object
$\wt G\in\ct^b_c$ we have that $\wt G$ belongs
to $\ct_H^{\leq n}$ for some $n$, hence
$\genul {\wt G}{}{0}\subset\ct_H^{\leq n}$,
and $\big[\genul {\wt G}{}{0}\big]^\perp\supset\ct_H^{\geq n+1}$.
Thus for any $\wt G\in\ct^b_c$ we have
$\ct^b_c\cap\big[\genul {\wt G}{}{0}\big]^\perp\supset\ct^b_c\cap\ct_H^{\geq n+1}$;
in other words the inequality
$\ct^b_c\cap\ct_H^{\geq0}
\preceq
\ct^b_c\cap\big[\genul {\wt G}{}{0}\big]^\perp$
is cheap.
\eprf

\lem{L1.95}
Let $\ct$ be a coherent triangulated category.
With $\cq(\ct^b_c)$ a minimal
$\ct^b_c\cap\big[\genul {G}{}{0}\big]^\perp$,
as in 
Definition~\ref{D0.13} and Lemma~\ref{L1.99}, the subcategory
$^\perp\cq(\ct^b_c)\cap\ct^b_c$ is equivalent, in the partial order, to
$\ct^b_c\cap\ct^{\leq0}$, where $\big(\ct^{\leq0},\ct^{\geq0}\big)$
is a \tstr\ in the preferred equivalence class.
\elem

\prf
By Lemma~\ref{L1.99} the subcategory $\cq(\ct^b_c)$
is equivalent to $\ct^b_c\cap\ct^{\geq0}$, hence
$^\perp\cq(\ct^b_c)\cap\ct^b_c$ is equivalent to
$^\perp\big[\ct^b_c\cap\ct^{\geq0}\big]\cap\ct^b_c$.
Clearly $\ct^b_c\cap\ct^{\leq-1}$ is a subcategory of
$^\perp\big[\ct^b_c\cap\ct^{\geq0}\big]\cap\ct^b_c$, that
is $\ct^b_c\cap\ct^{\leq0}\preceq{^\perp\cq(\ct^b_c)}\cap\ct^b_c$. We need to prove
the reverse inequality in the partial order.

Let $N>0$ be the integer of Definition~\ref{D29.1}, that is
any object $X\in\ct^-_c$ admits a triangle $A\la X\la B$ with
$A\in\ct^-_c\cap\ct^{\leq N}$ and $B\in\ct^b_c\cap\ct^{\geq0}$. I assert
that, with this choice of $N$, we have
$^\perp\big[\ct^b_c\cap\ct^{\geq0}\big]\cap\ct^b_c\subset\ct^b_c\cap\ct^{\leq N}$.
To prove this take $X\in{^\perp\big[\ct^b_c\cap\ct^{\geq0}\big]}\cap\ct^b_c$ and
form the triangle $A\la X\la B$ above. As $X$ belongs
to $^\perp\big[\ct^b_c\cap\ct^{\geq0}\big]\cap\ct^b_c$ and $B$ belongs to
$\ct^b_c\cap\ct^{\geq0}$ the map $X\la B$ must vanish, making
$X$ an direct summand of $A\in\ct^b_c\cap\ct^{\leq N}$.
\eprf

We summarizing the lemmas in this section in
the following Proposition. Note: in the
statement of Proposition we will work
\emph{entirely in the category}
$\ct^b_c$; thus all perpendiculars are understood
to be taken inside $\ct^b_c$.

\pro{P1.105}
Let $\ct$ be a coherent triangulated
category.
Then with the partial order of Definition~\ref{D0.13} 
there is (up to equivalence) a
unique minimal subcategory $\cq(\ct^b_c)$ among the $\big[\genul
G{}{0}\big]^\perp\subset\ct^b_c$. 
If we define $\cl_i$ by the fomula 
$\cl_i\op={^\perp\big[\T^i
\cq(\ct^b_c)\big]}$, 
then the subcategories $\{\cl_i,\,i\in\nn\}$
form a good metric equivalent to the $\{\cm\op_i\}$  defined on
$\big[\ct^b_c\big]\op$
in Example~\ref{E20.3}(ii).
\epro

\bibliographystyle{amsplain}
\bibliography{stan}

\end{document}